\documentclass{article}
\usepackage[utf8]{inputenc}
\usepackage{amsmath,amssymb}
\usepackage[ruled,commentsnumbered,linesnumbered,vlined]{algorithm2e}
\usepackage{xcolor}
\usepackage{indentfirst}
\usepackage{comment}
\usepackage{longtable}
\usepackage{lscape}
\usepackage{rotating}
\usepackage[title]{appendix}
\usepackage{amssymb}

\usepackage{multirow}

\usepackage{pdflscape}

\usepackage{lscape}
\usepackage{longtable}
\usepackage{amsthm}
\usepackage{subfigure}
\usepackage{picture}
\usepackage{graphicx}
\usepackage{fullpage}
\usepackage{float}

\usepackage{apacite}

\DeclareMathAlphabet\mathbfcal{OMS}{cmsy}{b}{n}

\title{A hybrid heuristic for capacitated three-level lot-sizing and replenishment problems with a distribution structure}

\author{ 
Jesus O. Cunha {\thanks{Departamento de Estatística e Matemática Aplicada, Universidade Federal do Ceará, Fortaleza, CE 60440-900, Brazil.  ({\tt jesus.ossian@dema.ufc.br})}}
    \and
    Geraldo R. Mateus {\thanks{Departamento de Ciência da Computação, Universidade Federal de Minas Gerais, Av. Antônio Carlos 6627, Belo Horizonte, MG 31270-010, Brazil.  ({\tt mateus@dcc.ufmg.br})}}
    \and
    Rafael A. Melo {\thanks{Institute of Computing, Universidade Federal da Bahia, Salvador, BA 40170-115, Brazil. ({\tt rafael.melo@ufba.br})} \thanks{Part of this work was developed while the author was a visiting scholar at the Departmento de Ciência da Computação, Universidade Federal de Minas Gerais, Brazil.}}
}
\begin{document}

\maketitle

\begin{abstract}

In this paper, we consider the capacitated three-level lot-sizing and replenishment problem with a distribution structure (3LSPD-C), recently proposed in the literature.
In 3LSPD-C, a single production plant delivers items to the warehouses from where they are distributed to their corresponding retailers. There is a capacity on the total amount the plant can produce in each period, whereas there are no capacities on transportation. The goal of this optimization problem consists in determining an integrated three-echelon production and distribution plan minimizing the total cost, which is composed of fixed costs for production and transportation setups as well as variable inventory holding costs. 
Additionally, we consider a generalization of the problem which also establishes storage capacities (or inventory bounds) on the warehouses and/or retailers, given the importance of such characteristics in practical industrial and commercial environments.
Such extension is denoted generalized capacitated three-level lot-sizing and replenishment problem with a distribution structure (G3LSPD-C).
We propose a hybrid mixed integer programming (MIP) heuristic consisting of a relax-and-fix approach to generate initial feasible solutions and a fix-and-optimize improvement procedure grounded on varying-size neighborhoods to obtain high-quality solutions. 
Computational experiments are performed to analyze the potential cost reductions achieved using the new hybrid heuristic when compared with a state-of-the-art MIP formulation.
The results show that the proposed hybrid heuristic can match or improve the solution quality obtained by the MIP formulation for the majority of the 3LSPD-C instances within a specified time limit.
Additionally, such superior behavior remains valid when the approaches are applied to the more general G3LSPD-C.
Furthermore, we investigate the economic impacts of the storage capacities and how they affect the performance of our newly proposed hybrid heuristic.
\\


\end{abstract}

\section{Introduction}
\label{sec:introduction}

\begin{figure}[H]
\centering
\includegraphics[scale=0.25]{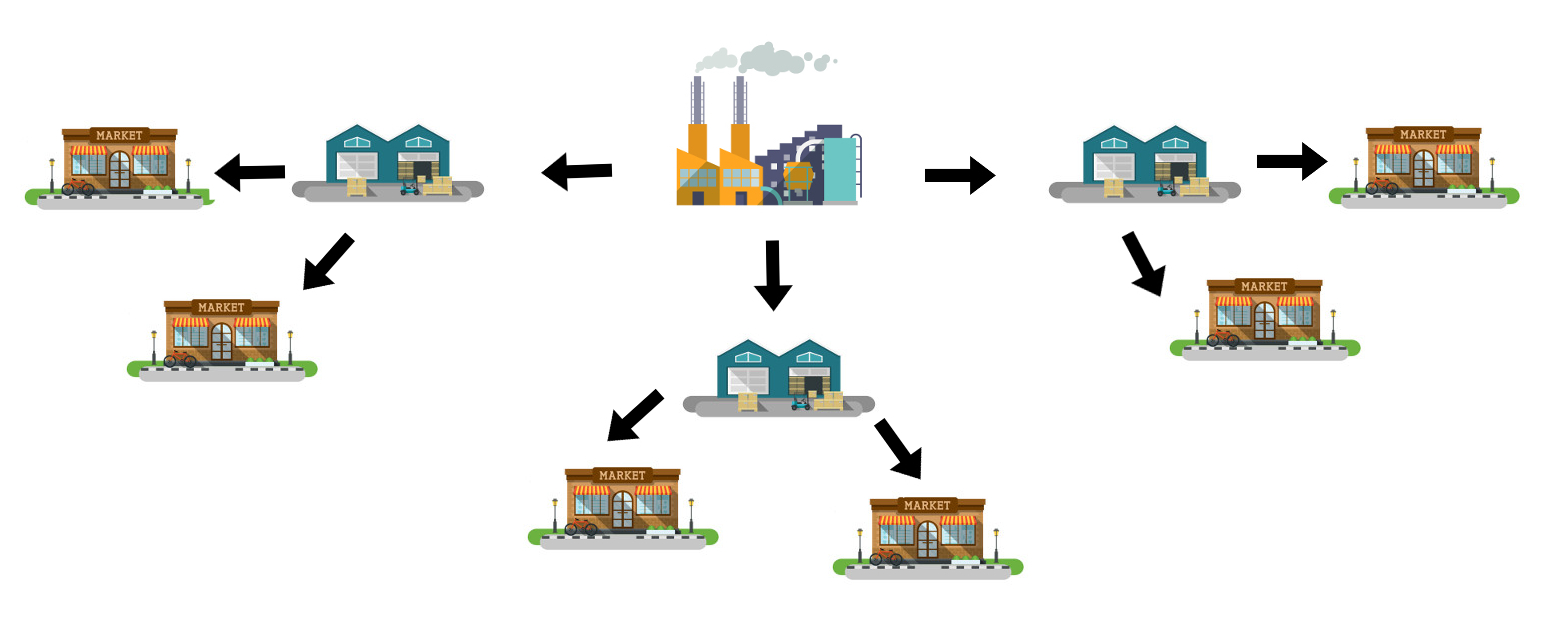}
\caption{Example of a three-level supply chain with a distribution structure.}
\label{fig:threelevelexample}
\end{figure}

Supply chain optimization plays a key role in the competitiveness of enterprises nowadays.
An important matter that can considerably increase effectiveness and reduce costs is an integrated optimization of the production and transportation plannings, a common need in companies with decentralized markets.
In this direction, we consider the capacitated three-level lot-sizing and replenishment problem with a distribution structure (3LSPD-C), introduced in \citeA{GruBazCorJan19}. In this three-echelon supply chain optimization problem, a single production plant produces items to replenish multiple warehouses, from where they are sent to fulfill the time-varying deterministic demands of several retailers over a finite planning horizon.
Additionally, there is a restriction on the quantity of items the plant can produce in each period.
The  goal consists in  determining  an  integrated  production, inventory, and distribution  plan which minimizes the total costs, consisting of fixed production and transportation setups together with per unit inventory holding costs.
Figure~\ref{fig:threelevelexample} exemplifies the layout of a three-level lot-sizing and replenishment problem with a distribution structure.
Notice that, in practical industrial and commercial settings, it is very natural that the storage facilities have capacity restrictions. In this regard, we also introduce the generalized capacitated three-level lot-sizing and replenishment problem with a distribution structure (G3LSPD-C), in which capacities on storage can occur at the warehouses and/or retailers.

Several problems involving integrated production and transportation have been studied in the literature.
Among these, we can highlight multi-level lot-sizing problems \cite{EppMar87,SlaBenDolMas20}, one-warehouse multi-retailer problems~\cite{SolSur12,CunMel16owmr} and some of its generalizations~\cite{Par05,MelWol12,LiHai19}, problems in the context of production and transportation by third-party logistic companies~\cite{SteZha07,MelWol10}, and supply chain optimization problems involving multiple plants in which items can be transferred between the plants~\cite{TanJha17,CarNas18,CunKraMel21}.

Recently, the optimization of three-echelon supply chains has received considerable attention.
\citeA{CarTre14} propose an integer programming formulation for a multi-product three-level supply chain problem.
\citeA{ZhaSon18} present an integer programming approach to tackle a practical three-level problem arising in a case study at Danone Waters China Division.
\citeA{GruBazCorJan19} introduce the uncapacitated and capacitated three-level lot-sizing and replenishment problem with a distribution structure (3LSPD-U and 3LSPD-C, respectively). The authors propose and compare, both theoretically and computationally, several formulations for both variants.
\citeA{GruCorJan20} describe a method based on Benders decomposition for a stochastic three-level lot-sizing and replenishment problem with a distribution structure.
\citeA{CunMel21} tackle 3LSPD-U and propose valid inequalities, preprocessing, and a multi-start randomized bottom-up dynamic programming-based heuristic. The authors managed to solve to optimality all the benchmark instances for the problem by combining the proposed heuristic with a preprocessed multicommodity formulation.

Production planning problems with storage capacities have appeared in several contexts.
It is noticeable, however, that most of the works related to storage capacities are very recent when compared to those considering production capacities~\cite{BitYan82,TriThoMcC89, BusSahHelTem10, VinDuhRenTch20,TavMin20}.
\citeA{AkbPenRap15} study capacitated lot-sizing problems with storage capacities. The authors propose polynomial-time algorithms for the single-item variant of the problem and show NP-hardness results for the multi-item variation.
\citeA{AkbPenRap15b} establish NP-hardness results for the multi-item uncapacitated lot-sizing with storage capacities under different scenarios of costs and capacities.
\citeA{BraAbsDauzKed15} propose integer programming formulations and a Lagrangian relaxation-based heuristic for an uncapacitated two-level lot-sizing problem with storage capacities.
\citeA{MelRib17} describe mixed integer programming (MIP) formulations and MIP heuristics for the multi-item uncapacitated lot-sizing problem with storage capacities.
\citeA{JinMu20} tackle a lot-sizing problem with perishable items, product substitution, and storage capacities. The authors provide a dynamic programming algorithm for the problem and analyze the effects of the different characteristics on the costs of the solutions.



\subsection{Main contributions and organization}

The main contribution of this paper is an effective hybrid heuristic for the capacitated three-level lot-sizing and replenishment problem with a distribution structure (3LSPD-C), which to the best of our knowledge is the first heuristic for the problem. 
The proposed hybrid heuristic combines a relax-and-fix approach with a fix-and-optimize improvement procedure which relies on varying-size neighborhoods to improve the local search mechanism and overcome local optima solutions.
Computational experiments evidence the effectiveness of the heuristic, as it outperforms a state-of-the-art formulation available for the problem.
The heuristic is flexible and can be easily extended to tackle more general variants of the problem, such as the generalized capacitated three-level lot-sizing and replenishment problem with a distribution structure (G3LSPD-C) that we introduce in our work.
G3LSPD-C considers storage capacities on the warehouses and/or retailers, which are very relevant characteristics as they often arise in commercial and industrial settings.
Last but not least, we analyze the economic impacts (in terms of the increase in the solutions' costs) of storage capacities on the warehouses and retailers.

The remainder of this paper is organized as follows. 
Section~\ref{sec:problemdefinition} formally defines the considered capacitated three-level lot-sizing and replenishment problems with a distribution structure.
Section~\ref{sec:proposedapproaches} describes the proposed hybrid heuristic.
Section~\ref{sec:experiments} summarizes the performed computational experiments.
Section \ref{sec:concludingremarks} concludes the paper with final remarks.

\section{Capacitated three-level lot-sizing and replenishment problems with a distribution structure}
\label{sec:problemdefinition}

In this section, we define two capacitated three-level lot-sizing and replenishment problems with a distribution structure. 
Section~\ref{sec:3lspdc} formalizes the capacitated three-level lot-sizing and replenishment problem with a distribution structure (3LSPD-C) studied in \citeA{GruBazCorJan19}. An available standard mixed integer programming formulation for the problem is also described.
Section~\ref{sec:g3lspdc} introduces the generalized capacitated three-level lot-sizing and replenishment problem with a distribution structure (G3LSPD-C), which extends 3LSPD-C by considering storage capacities on the warehouses and/or retailers.

\subsection{The capacitated three-level lot-sizing and replenishment problem with a distribution structure}\label{sec:3lspdc}

3LSPD-C can be formally stated as follows.
There is a set $F = P \cup W \cup R$ of facilities composed of a single plant, $P=\{p\}$, a set $W$ of warehouses, and a set $R$ of retailers.
Each retailer $r\in R$ has a time-varying deterministic demand $d^r_t$ for a single item for each period $t\in T$.
There is a capacity $C^p_t$ on the amount that can be produced at the production plant in period $t\in T$.
Each warehouse $w\in W$ attends a set of retailers $\delta(w)$ and each retailer $r\in R$ has a unique associated warehouse $\delta_w(r)$ from which it receives its items.
Each facility $i\in F$ has corresponding setup costs $sc^i_t$ and holding costs $hc^i_t$ for each period $t\in T$. The problem consists in encountering a feasible production/distribution plan which attends to all the demands without backlogging while minimizing the total cost.
Denote the demand in period $t \in T$ of the plant as $d^{p}_t = \sum_{r \in R} d^r_t$ and that of a warehouse $w\in W$ as $d^{w}_t = \sum_{r \in \delta(w)} d^r_t$. In addition, let the cumulative demand of facility $i\in F$ for periods $1\leq t \leq k \leq |T|$ be $d^i_{tk}=\sum_{l=t}^k d^i_l$. 

Consider the decision variable $x^i_t$ to be the amount produced in the production plant in period $t\in T$ if $i=p$ and to be equal to the amount transported to facility $i$ from its predecessor in period $t$ if $i\in W\cup R$. Besides, let the variable $s^i_t$ represent the amount of inventory in facility $i \in F$ in period $t\in T$. Additionally, define variable $y^i_t$ for $i\in F$ and $t\in T$ to be equal to one if $x^i_t>0$ and to be equal to zero otherwise. 3LSPD-C can be cast as the following standard MIP formulation~\cite{GruBazCorJan19}:
\begin{align}
(\textrm{3LSPD-C}) \qquad  & \  \min \ \ \  \sum_{t \in T} \left( \sum_{i\in F} sc^i_t y^i_t + \sum_{i\in F} hc^i_t s^i_t \right)  \label{std-obj} & \\
 \qquad & s^i_{t-1} + x^i_t = \sum_{j \in \delta(i)} x^j_t + s^i_t, \qquad  \textrm{for} \ i \in P\cup W, \ t\in T, \label{std-1} \\
&  s^r_{t-1} + x^r_t = d^r_t + s^r_t, \qquad  \textrm{for} \ r \in R, \ t\in T, \label{std-2} \\
&  x^i_t \leq d^i_{t|T|} y^i_t, \qquad  \textrm{for} \ i \in W\cup R, \ t\in T, \label{std-3} \\
&  x^p_t \leq \min\{C^p_t,d^p_{t|T|}\} y^p_t, \qquad  \textrm{for} \ \ t\in T, \label{std-3b} \\
&  x^i_{t}, \ s^i_t \geq 0, \qquad  \textrm{for} \ i \in F, \ t\in T, \label{std-4}\\
&  y^i_{t} \in  \{0,1\}, \qquad  \textrm{for} \ i \in F, \ t\in T. \label{std-5}
\end{align}
The objective function~\eqref{std-obj} minimizes the total cost, composed of setup and storage costs.
Constraints~\eqref{std-1}~and~\eqref{std-2} define inventory balance constraints for each of the facilities.
Constraints~\eqref{std-3}~and~\eqref{std-3b} are setup enforcing constraints, with~\eqref{std-3b} also guaranteeing the plant's production capacities are not violated. 
Constraints~\eqref{std-4} and \eqref{std-5} define, respectively, nonnegativity and integrality requirements on the variables.

\subsection{The generalized capacitated three-level lot-sizing and replenishment problem with a distribution structure}\label{sec:g3lspdc}

G3LSPD-C generalizes 3LSPD-C by considering the fact that the warehouses and/or retailers can have a limit on the amount that can be held in storage in any given period.
Define $\hat{C}^i_t$ to be the storage capacity at facility $i \in W \cup R$ in period $t\in T$.

G3LSPD-C can be formulated as the following MIP:
\begin{align}
(\textrm{G3LSPD-C}) \qquad  &   \min \ \ \  \sum_{t \in T} \left( \sum_{i\in F} sc^i_t y^i_t + \sum_{i\in F} hc^i_t s^i_t \right)  \label{estd-obj} & \\
 & \eqref{std-1}-\eqref{std-5}, \nonumber \\
&  s^i_{t} \leq \min\{\hat{C}^i_t, d^i_{t|T|} \}, \qquad  \textrm{for} \ i \in W\cup R, \ t\in T.  \label{estd-1} 
\end{align}
The objective function~\eqref{estd-obj} is the same as \eqref{std-obj}.
Constraints~\eqref{estd-1} ensure the total inventory at each facility $i \in W\cup R$ at the end of each period does not exceed the storage capacity.

\section{The hybrid heuristic}
\label{sec:proposedapproaches}

In this section, we detail the hybrid heuristic combining a relax-and-fix approach with a fix-and-optimize improvement procedure, both based on rolling horizon schemes.
This hybrid heuristic is motivated by the fact that MIP heuristics using these concepts have been successfully applied to many supply chain optimization problems~\cite{HelSah10,Che15,Toletal15,MelRib17,CuKrMe19}.
The proposed heuristic aims at providing good quality solutions while allowing flexibility to be easily adapted to variants of the problem.
Section~\ref{sec:relaxandfix} describes the relax-and-fix approach to obtain an initial feasible solution. Section~\ref{sec:fixandoptimize} explains the fix-and-optimize improvement procedure. The complete hybrid heuristic is presented in Section~\ref{sec:completehybrid}.

In what follows, we say that a variable belongs to an interval $[\alpha,\beta]$ if it corresponds to a period $t \in [\alpha,\beta]$.

\subsection{Relax-and-fix approach}
\label{sec:relaxandfix}

The proposed relax-and-fix approach consists of progressively building a solution by solving a series of subproblems in which a subset of the variables is fixed and another subset has its integrality requirements relaxed. It follows a rolling horizon mechanism, which works as follows. Consider integer values $\alpha^{RF}$ and $\beta^{RF}$ to determine three (possibly empty) intervals, corresponding to how the integer variables are treated when solving the mixed integer programs. The interval  $[1,\alpha^{RF}-1]$ corresponds to variables to be fixed, the interval $[\alpha^{RF},\beta^{RF}]$ determines the variables to be optimized with the integrality preserved, and the interval $[\beta^{RF}+1,|T|]$ defines the interval for which the integrality requirements on the variables are relaxed. The approach initiates with $[\alpha^{RF},\beta^{RF}]$ starting at the beginning of the planning horizon, and it progressively shifts towards the end of the horizon.
The relax-and-fix approach is described in Algorithm~\ref{alg:relaxfix}. 
It takes as inputs a mixed integer programming formulation FORM, an integer $k^{RF}$ defining the size of the horizon of the subproblems, an integer $k_{f}^{RF}$
representing the size of the interval to be fixed after solving a subproblem, and the maximum allowed time $maxt^{RF}$ for the approach.

\begin{algorithm}[H]
\caption{Relax-and-Fix (FORM, $k^{RF}$, $k_{f}^{RF}$, $maxt^{RF}$)}
\label{alg:relaxfix}
    $\Phi \leftarrow \varnothing$\; \label{RF-1}
     $\alpha^{RF} \leftarrow 1,\ \beta^{RF} \leftarrow \alpha^{RF} + k^{RF} - 1$\; \label{RF-3}
     $maxt^{RF}_{sub} \leftarrow \frac{maxt^{RF}}{\left\lceil |T|/k_{fix}^{RF}\right\rceil}$\; \label{RF:tlimitdef}
    \While{complete horizon was not yet treated}{ \label{RF:whilebeg}
        $(\hat{s},\hat{x},\hat{y}) \leftarrow $ Solve MIP$_{RF}$(FORM, $\alpha^{RF}$, $\beta^{RF}$, $\Phi$, $maxt^{RF}_{sub}$)\; \label{RF:solvemip}
         Update $\Phi$ with the fixings for the variables in the interval $[\alpha^{RF}, \alpha^{RF} + k^{RF}-1]$ according to the values of $\hat{y}$ \; \label{RF-7}
        $\alpha^{RF} \leftarrow \alpha^{RF} + k^{RF},\ \beta^{RF} \leftarrow \min\{\alpha^{RF}+k_{fix}^{RF}-1,|T|\}$\; \label{RF:updateinterval}
    } 
    \Return{solution $S = (\hat{s},\hat{x},\hat{y})$ determined by the fixings in $\Phi$, $elapt^{RF}$}\; \label{RF:return}
\end{algorithm}

The set of fixings imposed on the integer $y$ variables, $\Phi$, which defines a partial solution being built, is initialized as empty in line \ref{RF-1}. 
The values of $\alpha^{RF}$ and $\beta^{RF}$ are initialized in line \ref{RF-3}. The time limit for each subproblem, $maxt^{RF}_{sub}$, is determined in line \ref{RF:tlimitdef}.
The loop of lines~\ref{RF:whilebeg}-\ref{RF:updateinterval} is executed while the complete horizon was not yet treated, i.e., a complete integer feasible solution was not yet obtained.
Line~\ref{RF:solvemip} solves MIP$_{RF}$(FORM, $\alpha^{RF}$, $\beta^{RF}$, $\Phi$, $maxt^{RF}_{sub}$), in which the integer variables related to the periods in the interval $[1,\alpha^{RF} -1]$ are fixed according to the values determined in the previous rounds, which are defined by $\Phi$, those related to the periods in the interval $[\alpha^{RF},\beta^{RF}]$ must preserve the integrality, and the integrality constraints are relaxed for those related to the interval $[\beta^{RF} + 1, |T|]$.
In line~\ref{RF-7}, $\Phi$ is updated according to the return of the call to MIP$_{RF}$ in line~\ref{RF:solvemip}.
We consider two possible strategies for determining the fixing of the variables. The first strategy (S1) consists of adding to $S$ only the fixings of the variables which assume value 1, which implies that the MIP subproblems to be solved do not necessarily have all the variables in the interval $[1,\alpha^{RF}-1]$ fixed. The second one (S2) consists of inserting into $S$ the fixings corresponding to all the $y$ variables in the interval.
The values of $\alpha^{RF}$ and $\beta^{RF}$ are updated in line \ref{RF:updateinterval}.
The complete solution $S$ determined according to the  fixings in $\Phi$ and the elapsed time to run the approach are returned in line \ref{RF:return}. Note that $S$ corresponds to the solution obtained in the last call to MIP$_{RF}$ in line \ref{RF:solvemip}.

Observe that strategy S1 ensures a feasible solution can be encountered at the end of the execution of Algorithm~\ref{alg:relaxfix} whenever the instance is feasible and enough time is available. This is a consequence of the fact that there always exists a feasible solution to each call to MIP$_{RF}$. On the other hand, note that strategy S2 may lead to infeasible subproblems and, thus, does not ensure obtaining a feasible solution at the end of Algorithm~\ref{alg:relaxfix}.
We remark, though, that such infeasibilities are more likely to arise in situations with very strict capacities and/or very high demands concentrated in a few periods.

\subsection{Fix-and-optimize improvement procedure}
\label{sec:fixandoptimize}

We now propose a multi-round fix-and-optimize improvement procedure with varying-size neighborhoods.
The main intuition behind the varying-size neighborhoods is to improve the search mechanism using the concepts of variable neighborhood search~\cite{HanMlaBriMor19} and to attempt turning the approach less parameter-sensitive.
It follows a rolling-horizon approach which can restart from the beginning of the horizon, and its main idea works as follows. 
Consider values $\alpha^{FO}$ and $\beta^{FO}$ defining an interval $[\alpha^{FO},\beta^{FO}]$ determining variables which will be optimized. Additionally, define intervals $[1,\alpha^{FO}-1]$ and $[\beta^{FO}+1,|T|]$ corresponding to variables which will be fixed according to an available input solution $S$.
Similar to the relax-and-fix approach, the procedure initiates with $[\alpha^{FO},\beta^{FO}]$ starting at the beginning of the planning horizon, and it progressively shifts towards the end of the horizon. The difference lies in the fact that $[\alpha^{FO},\beta^{FO}]$ can be restarted at the beginning of the planning horizon and that the size of the intervals can be dynamically changed during the execution.
The procedure is described in Algorithm~\ref{alg:fixopt}. It receives as inputs a MIP formulation FORM, an initial feasible solution $S$, the minimum size for the optimization horizon $k^{FO}_{h,\min}$, the minimum size for the fixing interval $k^{FO}_{f,\min}$, the length of the increase in the optimization horizon $\Delta^{FO}_h$, the length of the increase in the fixing horizon $\Delta^{FO}_f$, the minimum number of rounds to be performed $rnd^{FO}_{\min}$, and a time limit $maxt^{FO}$.
The algorithm assumes that $k^{FO}_{h,\min} \geq k^{FO}_{f,\min}$ and $\Delta^{FO}_h \geq \Delta^{FO}_f$.

\begin{algorithm}[H]
\caption{Fix-and-Optimize (FORM, ${S}$, $k^{FO}_{h,\min}$, $k^{FO}_{f,\min}$, $\Delta^{FO}_h$, $\Delta^{FO}_f$, $rnd^{FO}_{\min}$, $maxt^{FO}$)}
\label{alg:fixopt}
     $S^*,S'\leftarrow S$\;\label{FO-1}
     $k^{FO} \leftarrow k^{FO}_{h,\min}$,\ \ $k^{FO}_{f} \leftarrow k^{FO}_{f,\min}$ \;\label{FO-5}
     $maxt^{FO}_{sub} \leftarrow  \dfrac{maxt^{FO}}{rnd^{FO}_{\min} \times \left\lceil \frac{|T|}{k^{FO}_{\min}}\right\rceil } $\;\label{FO-3}
    \Repeat{stopping criterion is met}{ \label{FO:repeatbeg}
         
          $\alpha^{FO} \leftarrow 1,\ \ \beta^{FO} \leftarrow \min \{ \alpha^{FO} + k^{FO} - 1,\ |T| \}$\; \label{FO-6}
        \While{$elapt^{FO} < maxt^{FO}$ and complete horizon was not yet treated}{ \label{FO:whilebeg}
             $S \leftarrow$ MIP$_{FO}$(FORM, $\alpha^{FO}$, $\beta^{FO}$, $S^{*}$, $maxt^{FO}_{sub}$)\; \label{FO:solvemip}
            \If{$z(S) < z(S^{*})$}{\label{FO-ifimprove}
                 $S^* \leftarrow S$\; \label{FO:updatesol}
            }        
             $\alpha^{FO} \leftarrow \min \{\alpha^{FO} + k^{FO}_{f}, |T|\},\ \ \beta^{FO} \leftarrow \min\{\alpha^{FO}+k^{FO}-1,|T|\}$\; \label{FO:updateinterval}
             \If{$elapt^{FO}< maxt^{FO}$ and $k^{FO} < |T|$ }{\label{FO:ifupdatemaxtsub}
                $maxt^{FO}_{sub} \leftarrow \min\{maxt^{FO}_{sub},maxt^{FO}- elapt^{FO}\}$\; \label{FO:updatemaxtsub}
            }
        } \label{FO:whileend}
        \uIf{$z(S^*) < z(S')$}{\label{FO-7}
             $S'\leftarrow S^*$\; \label{FO-8}
        }\uElseIf{$k^{FO} < |T|$}{\label{FO-9}
             $k^{FO} \leftarrow k^{FO} + \Delta^{FO}_h$\; \label{FO-10}
             $k^{FO}_f \leftarrow k^{FO}_f + \Delta^{FO}_f$\; \label{FO-10b}
             \uIf{$k^{FO}\geq |T|$}{\label{FO:ifmaxtime}
                $maxt^{FO}_{sub} \leftarrow maxt^{FO}- elapt^{FO}$\;\label{FO:setmaxtime}
             }
        }\Else{   
             stop\; \label{FO-11}
        }
    } \label{FO:repeatend}
    \Return{best obtained solution $S^*$}\;\label{FO:return}
\end{algorithm}

Initially, line~\ref{FO-1} defines the input solution $S$ as both the best and previously best-known solutions ($S^*$ and $S'$, respectively).
Line~\ref{FO-5} initializes the sizes of the optimization and fixing intervals ($k^{FO}$ and $k^{FO}_{f}$, correspondingly).
Line~\ref{FO-3} sets the maximum time for solving each subproblem, which ensures that at least $rnd^{FO}_{\min}$ rounds are performed.
Each iteration of the repeat loop (lines~\ref{FO:repeatbeg}-\ref{FO:repeatend}) represents one round of the fix-and-optimize procedure. 
In each round, the algorithm goes through the whole planning horizon in order to solve the corresponding subproblems.
Firstly, the values $\alpha^{FO}$ and $\beta^{FO}$, defining the interval related to the subproblem to be optimized, are initialized in line~\ref{FO-6}.
The while loop (lines~\ref{FO:whilebeg}-\ref{FO:whileend}) represents the iterations performed in each round. The subproblem corresponding to the interval $[\alpha^{FO},\beta^{FO}]$ and solution $S^*$ is solved in line~\ref{FO:solvemip}, and in case the current best solution is improved, it is updated in line~\ref{FO:updatesol}. After that, the interval $[\alpha^{FO},\beta^{FO}]$ is updated in line~\ref{FO:updateinterval}.
Lines~\ref{FO:ifupdatemaxtsub}-\ref{FO:updatemaxtsub} update the allowed time for solving the subproblems, whenever appropriate.
Next, if the current best solution improves over the previous best, the latter is updated (lines~\ref{FO-7}-\ref{FO-8}). Otherwise, if the size of the fixing interval did not achieve the size of the horizon, the sizes of the optimization and fixing intervals are updated.
In case the size of the optimization interval reached that of the planning horizon, the remaining allowed time is defined as the time for solving the only subproblem to be tackled (lines~\ref{FO:ifmaxtime}-\ref{FO:setmaxtime}).
If none of the two previous conditions are satisfied, the approach has reached a stopping criterion.
The best-obtained solution $S^*$ is returned in line~\ref{FO:return}.

We remark that the choices of the control parameters are crucial for the good behavior of the proposed fix-and-optimize improvement procedure. Note that $k^{FO}_{h,\min}$ should define a reasonable size for the subproblems to be solved to optimality or near optimality in the first rounds. Besides, the size $k^{FO}_{f,\min}$ should allow flexibility in the reoptimization of latter periods in a given round. Furthermore, $rnd^{FO}_{min}$ should be determined so that it allows enough time for solving the subproblems.

\subsection{The hybrid heuristic}
\label{sec:completehybrid}

The hybrid heuristic consists of generating an initial feasible solution using the relax-and-fix approach, described in Section~\ref{sec:relaxandfix}, followed by a call to the fix-and-optimize improvement procedure, described in Section~\ref{sec:fixandoptimize}.
The heuristic is presented in Algorithm~\ref{alg:hybridheuristic} and it takes as inputs the parameters corresponding to the calls to Relax-and-Fix (Algorithm~\ref{alg:relaxfix}) and Fix-and-Optimize (Algorithm~\ref{alg:fixopt}), as well as a time limit $maxt$.

\begin{algorithm}[H]
\caption{Hybrid-Heuristic (FORM, $k^{RF}$,$k_{f}^{RF}$, $maxt^{RF}$,$k^{FO}_{h,\min}$, $k^{FO}_{f,\min}$, $\Delta^{FO}_h$, $\Delta^{FO}_f$, $rnd^{FO}_{\min}$, $maxt$)}
\label{alg:hybridheuristic}
    $S,elapt \leftarrow$ Relax-and-Fix (FORM, $k^{RF}$, $k_{f}^{RF}$, $maxt^{RF}$)\; \label{HH-1}
     $S^* \leftarrow$  Fix-and-Optimize (FORM, ${S}$, $k^{FO}_{h,\min}$, $k^{FO}_{f,\min}$, $\Delta^{FO}_h$, $\Delta^{FO}_f$, $rnd^{FO}_{\min}$, $maxt - elapt$)\; \label{HH-2}
    \Return{$S^*$}\; \label{HH:return}
\end{algorithm}

\section{Computational experiments}
\label{sec:experiments}

This section reports the computational experiments conducted to evaluate the performance of the proposed hybrid heuristic.
All computational experiments were carried out on a machine running under Ubuntu GNU/Linux, with an Intel(R) Core(TM) i5-3740 CPU @ 3.20GHz processor and 8Gb of RAM. The algorithms were coded in Julia v1.6.2, using JuMP v0.21.8. The formulations were solved using Gurobi 9.0.2 with the standard configurations, except the relative optimality tolerance gap, which was set to $10^{-6}$.

\subsection{Benchmark instances}\label{sec:instancias}

\subsubsection*{3LSPD-C}

The benchmark instances for 3LSPD-C were introduced in \citeA{GruBazCorJan19}.
All instances have a single plant. For each instance, the number of retailers $|R|$ lies in $\{50, 100, 200\}$, while that of warehouses belongs to $\{5, 10, 15, 20\}$. 
The demand for each period is randomly defined using a uniform distribution $U[5, 100]$.
Fixed setup costs are defined as follows: $sc^p$ is uniformly chosen from $U[30000, 45000]$, $sc^{w}$ from $U[1500, 4500]$, and $sc^r$ from $U[5, 100]$. The storage costs are defined as $hc^p = 0.25$, $hc^{\delta_w(r)}=0.5$, and $hc^r$ is uniformly chosen from $U[0.5,1]$.
Instances are classified as unbalanced (in which a few warehouses serve most of the retailers) or balanced (in which all the warehouses serve nearly the same amount of retailers).
Considering a capacity factor $C \in \{1.50, 1.75, 2.00\}$, the time-invariant plant capacity for every period $t\in T$ is defined as $C^p_t = \frac{C}{|T|}\sum_{i\in R}\sum_{t\in T}d^i_t$.
We only consider in our experiments the instances with time-invariant fixed setups and dynamic demands. For each different combination of these parameters, there are five instances, giving a total of 360 unbalanced and 360 balanced instances.

\subsubsection*{G3LSPD-C}

For G3LSPD-C, we considered adaptations of the instances for 3LSPD-U and 3LSPD-C~\cite{GruBazCorJan19}.
Namely, for each instance of 3LSPD-U and 3LSPD-C, we generated instances considering a storage capacity factor $C_s \in \{1.50, 1.75, 2.00\}$ to define time-invariant facility storage capacities 
$\hat{C}^i_t = \frac{C}{|T|}\sum_{t\in T}d^i_t$ for every period $t \in T$. 
The considered instances either have storage capacities on the warehouses ($i\in W$) or the retailers ($i\in R$), but never on both.
Thus, the benchmark set is composed of 2880 unbalanced instances and 2880 balanced ones.

\subsection{Considered approaches and settings}
\label{sec:approachesandsettings}

The following approaches were compared in our experiments:
\begin{itemize}
    \item the newly proposed hybrid heuristic (RFFO), presented in Section~\ref{sec:proposedapproaches};
    \item the echelon stock formulation (ES-LS), as described in \cite{GruBazCorJan19}.
\end{itemize}
The comparison of RFFO with ES-LS is because the latter was shown to be amongst the most effective approaches for 3LSPD-C in~\citeA{GruBazCorJan19}. ES-LS is presented in Appendix~\ref{sec:echelonstockformulations} for the sake of completeness.

Preliminary experiments were carried out to determine the parameters for RFFO. They took into consideration 144 randomly selected instances that are not used as part of the final experiments.
Firstly, we highlight that RFFO employed the echelon stock formulation ES-LS~\cite{GruBazCorJan19}.
As a general observation, the fix-and-optimize improvement procedure (FO) was the main ingredient of the heuristic, as it could considerably improve the results of the relax-and-fix approach (RF), even when initial solutions with lower quality were obtained.
The preliminary experiments allowed us to notice that the performance of RFFO was not sensitive to the tested parameters for RF ($k^{RF}$ and $k_{f}^{RF}$). Besides, it was also not sensitive to the initial parameters of FO ($k^{FO}_{h,\min}$ and $k^{FO}_{f,\min}$) when reasonably small values were used. The results showed that values of $rnd^{FO}_{\min} \in \{2,3\}$ achieved similar results, which were better than those using other smaller or larger values.
At the end, the parameters for the complete experiments were defined as follows: 
$k^{RF}=5$, $k_{f}^{RF}=3$, $k^{FO}_{h,\min}=5$, $k^{FO}_{f,\min}=3$, $\Delta^{FO}_h = 0$, $\Delta^{FO}_f=1$, and $rnd^{FO}_{\min}=2$. 
The time limits for RFFO were set as $maxt=600$ seconds and $maxt^{RF}= \lceil 0.10 \times maxt \rceil = 60$ seconds.

The formulations were executed with a time limit of 600 seconds to allow a direct comparison with RFFO regarding the quality of the achieved solutions. Besides, for 3LSPD-C they were also executed with a time limit of 3600 seconds to verify whether additional instances could be solved to optimality.

\subsection{Results for 3LSPD-C}
\label{sec:results3lspdc}

The results for 3LSPD-C are summarized in Tables \ref{tab:lab0_600s}-\ref{tab:lab1_3600s}, where Tables~\ref{tab:lab0_600s}-\ref{tab:lab1_600s} compare the results of RFFO with those of ES-LS running for 600 seconds, while Tables~\ref{tab:lab0_3600s}-\ref{tab:lab1_3600s} compare our results with those running ES-LS for 3600 seconds. 
In Tables~\ref{tab:lab0_600s}-\ref{tab:lab1_600s}, the first three columns identify the instance parameters (param), their corresponding values (value), and the number of instances with the specified characteristics (\#inst). 
The next three columns provide, for ES-LS within the time limit of 600 seconds, the average best solution (best$_{ESLS}$), the average open optimality gap in percent (gap$_{ESLS}$, calculated for each instance as $100 \times \frac{\textrm{best}_{ESLS}-\textrm{bbound}}{\textrm{best}_{ESLS}}$, where bbound denotes the best dual bound achieved by the solver at the end of the time limit), and the number of instances solved to optimality (\#opt). 
The next eight columns show
the average best solution achieved by RF with a time limit of 60 seconds (best$_{RF}$),
the average best solution obtained at the end of RFFO within 600 seconds (best$_{RFFO}$),
the average open optimality gap in percent obtained by RFFO (gap$_{RFFO}$, given for each instance by $100 \times \frac{\textrm{best}_{RFFO} - \textrm{bbound}}{\textrm{best}_{RFFO}}$), 
the average gain in percent achieved by FO over RF (impr$_{\frac{FO}{RF}}$, calculated for each instance as $100 \times \frac{\textrm{best}_{RF} - \textrm{best}_{RFFO}}{\textrm{best}_{RF}}$), 
the average gain in percent achieved by RFFO over ES-LS (impr$_{\frac{RFFO}{ESLS}}$, calculated for each instance as $100 \times \frac{\textrm{best}_{ESLS} - \textrm{best}_{RFFO}}{\textrm{best}_{ESLS}}$),
the number of instances for which RFFO achieved a solution at least as good as that of ES-LS ($\leq$best$_{ESLS}$), 
the number of instances for which RFFO obtained a solution which strictly improves over that encountered by ES-LS ($<$best$_{ESLS}$),
and the average number of rounds performed by FO (rnd$_{FO}$).
The best average solutions in each line are highlighted in bold.
Tables~\ref{tab:lab0_3600s}-\ref{tab:lab1_3600s} present a subset of these just described columns, with the difference that those corresponding to ES-LS report information for the formulation with a time limit of 3600 seconds.

Tables \ref{tab:lab0_600s}-\ref{tab:lab1_600s} show the robustness of RFFO, as it obtained lower average best values than ES-LS running with a time limit of 600 seconds for all the instance configurations. 
RFFO achieved overall average improvements over ES-LS of 3.07\% and 2.95\% for the unbalanced and balanced instances, respectively. One can observe that the improvements reached values as high as 7.33\% and 7.93\% for the instances with a larger number of retailers ($|R|=200$).
It should also be noticed that the gaps of RFFO are remarkably lower than those of 
ES-LS.
Additionally, the number of instances for which RFFO at least matched the best solutions of ES-LS is remarkably high, with an overall value of 303 and 285 out of 360 instances (respectively, 84.17\% and 79.17\%).
Besides, FO allowed considerable improvements over the solutions obtained by RF, especially for the larger instances with $|T|=30$ and those with $|R|=200$.
Furthermore, it can be noticed that the parameter that influenced the most the average number of rounds performed by FO for both unbalanced and balanced instances was the number of retailers $|R|$.
When we compare the results for the unbalanced instances (Table~\ref{tab:lab0_600s}) with those for the balanced ones (Table~\ref{tab:lab1_600s}), we can see that the behaviors of ES-LS and RFFO do not change much when we analyze the achieved average gaps. In this case, though, it is noteworthy that RFFO encountered a larger number of solutions at least matching those of ES-LS for the unbalanced instances, whereas a larger number of strictly improving solutions was obtained for the balanced ones.

\begin{figure}[H]
\centering
   \includegraphics[width=0.55\textwidth] {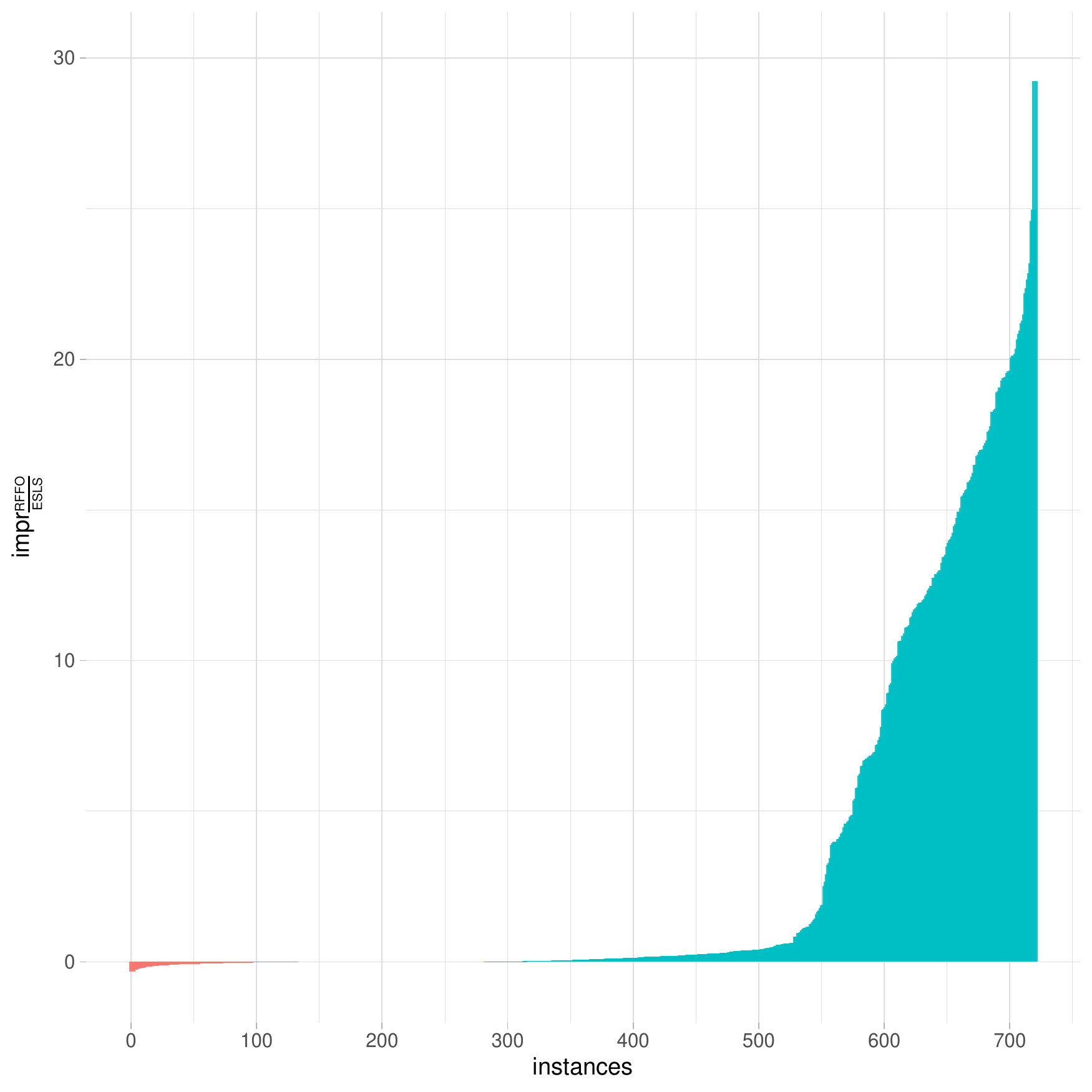}
\caption{Plot with the improvements of RFFO over ES-LS (impr$_{\frac{RFFO}{ESLS}}$) for all 3LSPD-C instances.}
\label{fig:improvements}
\end{figure}

The plot in Figure~\ref{fig:improvements} shows the improvements of RFFO over ES-LS for all the 720 instances (both balanced and unbalanced). It can be seen in the figure that the positive improvements of RFFO over ES-LS can be very significant as they reach close to 30\% in certain cases. On the other hand, the improvements of ES-LS over RFFO, determined by the negative values, are negligible.

Tables~\ref{tab:lab0_3600s}-\ref{tab:lab1_3600s} show that even running ES-LS for a much larger time (3600 seconds), RFFO still obtained good comparative numbers even though it was executed for much less time (600 seconds).
Note that RFFO achieved lower average best values for all configurations except for the smaller instances with $|T|=15$ and those with $R=50$.
Besides, remark that the average open optimality gaps achieved by RFFO are reasonably low (1.24\% and 1.30\% for the unbalanced and balanced instances, respectively), and these values reached a maximum average of 2.39\% for $|T|=30$ in the unbalanced instances and 2.18\% for $C=1.75$ in the balanced ones. 
We remark that the open optimality gaps for RFFO in these tables are smaller than those reported in Tables~\ref{tab:lab0_600s}-\ref{tab:lab1_600s} as the dual bounds obtained by ES-LS within 3600 seconds are better than those obtained with a time limit of 600 seconds.
When we compare the results for the unbalanced instances (Table~\ref{tab:lab0_3600s}) with those for the balanced ones (Table~\ref{tab:lab1_3600s}), it can be observed that the behavior of ES-LS does not change much for most of the instance configurations when we consider the achieved average gaps, except for the parameter $|W|$. This leads to similar behaviors of RFFO for unbalanced and balanced instances when compared to the solutions encountered by ES-LS, except for the improvements achieved for the different values of $|W|$, as that for $|W|=5$ shows a much higher average for the balanced instances.

\begin{landscape}

\begin{table}[H]
\centering
\footnotesize
\caption{Results for 3LSPD-C considering the unbalanced instances and ES-LS running for 600 seconds.}
\begin{tabular}{l c c|ccc|cccccccc} \hline
\multicolumn{3}{c|}{Instances}  &  \multicolumn{3}{c|}{ES-LS} &  \multicolumn{8}{c}{RFFO} \\
param & value & \#inst &    best$_{ESLS}$ &  gap$_{ESLS}$ &  \#opt &   best$_{RF}$  &   best$_{RFFO}$  &  gap$_{RFFO}$ &  impr$_{\frac{FO}{RF}}$ &  impr$_{\frac{RFFO}{ESLS}}$ &  $\leq$best$_{ESLS}$ & $<$best$_{ESLS}$ &  rnd$_{FO}$ \\ \hline
 ALL &      &  360 &  867441.64 & 4.38 &  105 &  888640.79 &  \textbf{823577.22} &      1.41 &        4.58 &            3.07 &      303 &     204 &    9.22 \\ [0.1cm]
   $|T|$ &   15 &  180 &  555507.71 & 0.28 &  105 &  557932.85 &  \textbf{555497.91} &      0.28 &        0.41 &            0.00 &      139 &      40 &   11.03 \\
     &   30 &  180 & 1179375.57 & 8.48 &    0 & 1219348.72 & \textbf{1091656.53} &      2.53 &        8.76 &            6.14 &      164 &     164 &    7.41 \\ [0.1cm]
   $|W|$ &    5 &   90 &  752824.44 & 4.28 &   32 &  763302.69 &  \textbf{718114.70} &      1.49 &        4.01 &            2.89 &       89 &      51 &   10.91 \\
     &   10 &   90 &  839057.40 & 4.58 &   32 &  850736.95 &  \textbf{794950.18} &      1.46 &        4.12 &            3.23 &       73 &      49 &    9.50 \\
     &   15 &   90 &  888552.48 & 4.43 &   22 &  912060.35 &  \textbf{841440.29} &      1.31 &        4.80 &            3.21 &       75 &      57 &    8.53 \\
     &   20 &   90 &  989332.25 & 4.24 &   19 & 1028463.18 &  \textbf{939803.72} &      1.36 &        5.40 &            2.96 &       66 &      47 &    7.93 \\ [0.1cm]
      $|R|$ & 50 &  120 &  717586.93 & 0.94 &   50 &  724468.88 &  \textbf{717099.07} &      0.89 &        0.79 &            0.06 &       96 &      47 &   11.55 \\
     & 100 &  120 &  834384.89 & 2.97 &   47 &  846863.10 &  \textbf{813536.51} &      1.18 &        3.00 &            1.83 &      104 &      65 &    9.54 \\
     & 200 &  120 & 1050353.10 & 9.23 &    8 & 1094590.39 &  \textbf{940096.09} &      2.15 &        9.96 &            7.33 &      103 &      92 &    6.57 \\[0.1cm]
   $C$ &   1.50 &  120 &  948724.51 & 4.49 &   20 &  956560.93 &  \textbf{891806.31} &      0.91 &        4.24 &            3.66 &       94 &      75 &    8.43 \\
     &  1.75 &  120 &  868440.03 & 5.04 &   40 &  890309.50 &  \textbf{825931.35} &      2.19 &        4.54 &            2.98 &      106 &      68 &    9.17 \\
     &  2.00 &  120 &  785160.39 & 3.61 &   45 &  819051.93 &  \textbf{752994.00} &      1.12 &        4.97 &            2.57 &      103 &      61 &   10.07 \\
\hline
\end{tabular}
\label{tab:lab0_600s}
\end{table}

\begin{table}[H]
\centering
\footnotesize
\caption{Results for 3LSPD-C considering the balanced instances and ES-LS running for 600 seconds.}
\begin{tabular}{l c c|ccc|cccccccc} \hline
\multicolumn{3}{c|}{Instances}  &  \multicolumn{3}{c|}{ES-LS} &  \multicolumn{8}{c}{RFFO} \\
 param & value  & \#inst &    best$_{ESLS}$ &  gap$_{ESLS}$ &  \#opt &   best$_{RF}$  &   best$_{RFFO}$  &  gap$_{RFFO}$ &  impr$_{\frac{FO}{RF}}$ &  impr$_{\frac{RFFO}{ESLS}}$ &  $\leq$best$_{ESLS}$ & $<$best$_{ESLS}$ &  rnd$_{FO}$ \\ \hline
 ALL &      &  360 &  916176.01 &  4.33 &   52 &  936703.77 &  \textbf{871189.82} &      1.47 &        4.38 &            2.95 &      285 &     237 &    9.03 \\[0.1cm]
   $|T|$ &   15 &  180 &  589491.91 &  0.76 &   52 &  591987.16 &  \textbf{589158.03} &      0.71 &        0.43 &            0.05 &      128 &      80 &   10.57 \\
     &   30 &  180 & 1242860.10 &  7.90 &    0 & 1281420.39 & \textbf{1153221.61} &      2.23 &        8.33 &            5.85 &      157 &     157 &    7.49 \\[0.1cm]
   $|W|$ &    5 &   90 &  809926.18 &  6.02 &   17 &  788994.27 &  \textbf{750806.41} &      1.79 &        3.32 &            4.40 &       81 &      66 &    9.81 \\
     &   10 &   90 &  893863.57 &  4.63 &   13 &  904856.00 &  \textbf{843910.50} &      1.59 &        4.30 &            3.14 &       77 &      66 &    9.36 \\
     &   15 &   90 &  944040.67 &  3.97 &   12 &  970783.24 &  \textbf{900289.37} &      1.39 &        4.55 &            2.64 &       70 &      60 &    8.86 \\
     &   20 &   90 & 1016873.60 &  2.72 &   10 & 1082181.58 &  \textbf{989753.00} &      1.12 &        5.33 &            1.63 &       57 &      45 &    8.10 \\[0.1cm]
   $|R|$ & 50 &  120 &  736511.85 &  0.90 &   43 &  743821.56 &  \textbf{735904.94} &      0.84 &        0.82 &            0.07 &       88 &      47 &   12.05 \\
     & 100 &  120 &  870414.42 &  1.96 &    9 &  898640.97 &  \textbf{861010.05} &      1.12 &        3.11 &            0.86 &       89 &      82 &    9.40 \\
     & 200 &  120 & 1141601.74 & 10.13 &    0 & 1167648.79 & \textbf{1016654.47} &      2.46 &        9.19 &            7.93 &      108 &     108 &    5.64 \\[0.1cm]
   $C$ &   1.50 &  120 &  992888.15 &  4.01 &    5 & 1006824.08 &  \textbf{939766.40} &      0.87 &        4.15 &            3.21 &       97 &      89 &    8.00 \\
     &  1.75 &  120 &  916194.64 &  4.99 &   19 &  937903.20 &  \textbf{873496.42} &      2.28 &        4.30 &            2.83 &       91 &      73 &    9.04 \\
     &  2.00 &  120 &  839445.23 &  3.99 &   28 &  865384.04 &  \textbf{800306.64} &      1.27 &        4.68 &            2.82 &       97 &      75 &   10.05 \\
\hline
\end{tabular}
\label{tab:lab1_600s}
\end{table}

\begin{table}[H]
\centering
\footnotesize
\caption{Results for 3LSPD-C considering the unbalanced instances and ES-LS running for 3600 seconds.}
\begin{tabular}{l c c|ccc|ccccc} \hline
\multicolumn{3}{c|}{Instances}  &  \multicolumn{3}{c|}{ES-LS} &  \multicolumn{5}{c}{RFFO} \\
 param & value  & \#inst &    best$_{ESLS}$ &  gap$_{ESLS}$ &  \#opt   &  best$_{RFFO}$  &  gap$_{RFFO}$ &   impr$_{\frac{RFFO}{ESLS}}$ &  $\leq$best$_{ESLS}$ & $<$best$_{ESLS}$ \\ \hline
 ALL &      &  360 &  836685.34 & 2.11 &  134 &  \textbf{823577.22} &      1.24 &            0.90 &      242 &     134 \\ [0.1cm]
   $|T|$ &   15 &  180 &  \textbf{555336.71} & 0.07 &  134 &  555497.91 &      0.09 &           -0.02 &      115 &       7 \\
     &   30 &  180 & 1118033.98 & 4.16 &    0 & \textbf{1091656.53} &      2.39 &            1.82 &      127 &     127 \\ [0.1cm]
   $|W|$ &    5 &   90 &  719033.14 & 1.45 &   33 &  \textbf{718114.70} &      1.37 &            0.09 &       78 &      37 \\
     &   10 &   90 &  814614.99 & 2.66 &   35 &  \textbf{794950.18} &      1.32 &            1.39 &       63 &      37 \\
     &   15 &   90 &  868882.24 & 2.93 &   35 &  \textbf{841440.29} &      1.15 &            1.83 &       55 &      33 \\
     &   20 &   90 &  944211.01 & 1.40 &   31 &  \textbf{939803.72} &      1.13 &            0.28 &       46 &      27 \\ [0.1cm]
      $|R|$ & 50 &  120 &  \textbf{716942.52} & 0.79 &   57 &  717099.07 &      0.80 &           -0.01 &       76 &      27 \\
     & 100 &  120 &  814314.07 & 1.18 &   54 &  \textbf{813536.51} &      1.12 &            0.07 &       87 &      46 \\
     & 200 &  120 &  978799.44 & 4.36 &   23 &  \textbf{940096.09} &      1.80 &            2.64 &       79 &      61 \\ [0.1cm]
   $C$ &   1.50 &  120 &  910793.00 & 1.88 &   34 &  \textbf{891806.31} &      0.73 &            1.18 &       69 &      48 \\
     &  1.75 &  120 &  834645.99 & 2.64 &   45 &  \textbf{825931.35} &      2.04 &            0.62 &       88 &      46 \\
     &  2.00 &  120 &  764617.04 & 1.81 &   55 &  \textbf{752994.00} &      0.94 &            0.90 &       85 &      40 \\ 
\hline
\end{tabular}
\label{tab:lab0_3600s}
\end{table}

\begin{table}[H]
\centering
\footnotesize
\caption{Results for 3LSPD-C considering the balanced instances and ES-LS running for 3600 seconds.}
\begin{tabular}{l c c|ccc|ccccc} \hline
\multicolumn{3}{c|}{Instances}  &  \multicolumn{3}{c|}{ES-LS} &  \multicolumn{5}{c}{RFFO} \\
 param & value & \#inst &    best$_{ESLS}$ &  gap$_{ESLS}$ &  \#opt & best$_{RFFO}$  &  gap$_{RFFO}$ &   impr$_{\frac{RFFO}{ESLS}}$ &  $\leq$best$_{ESLS}$ & $<$best$_{ESLS}$ \\ \hline
 ALL &      &  360 &  884173.08 & 2.18 &   93 &  \textbf{871189.82} &      1.30 &            0.92 &      219 &     155 \\ [0.1cm]
   $|T|$ &   15 &  180 &  \textbf{589045.88} & 0.42 &   93 &  589158.03 &      0.44 &           -0.02 &      101 &      37 \\
     &   30 &  180 & 1179300.27 & 3.93 &    0 & \textbf{1153221.61} &      2.15 &            1.85 &      118 &     118 \\ [0.1cm]
   $|W|$ &    5 &   90 &  798782.99 & 4.85 &   27 &  \textbf{750806.41} &      1.57 &            3.41 &       76 &      54 \\
     &   10 &   90 &  845040.99 & 1.46 &   24 &  \textbf{843910.50} &      1.39 &            0.08 &       60 &      42 \\
     &   15 &   90 &  902055.34 & 1.30 &   23 &  \textbf{900289.37} &      1.19 &            0.11 &       46 &      34 \\
     &   20 &   90 &  990812.99 & 1.09 &   19 &  \textbf{989753.00} &      1.03 &            0.06 &       37 &      25 \\ [0.1cm]
    $|R|$ & 50 &  120 &  \textbf{735806.22} & 0.76 &   55 &  735904.94 &      0.77 &           -0.01 &       69 &      24 \\
     & 100 &  120 &  861601.89 & 1.06 &   37 &  \textbf{861010.05} &      1.00 &            0.05 &       68 &      52 \\
     & 200 &  120 & 1055111.12 & 4.71 &    1 & \textbf{1016654.47} &      2.11 &            2.70 &       82 &      79 \\ [0.1cm]
   $C$ &   1.50 &  120 &  955254.73 & 1.70 &   18 &  \textbf{939766.40} &      0.73 &            0.99 &       71 &      56 \\
     &  1.75 &  120 &  887681.20 & 3.15 &   36 &  \textbf{873496.42} &      2.18 &            1.02 &       78 &      57 \\
     &  2.00 &  120 &  809583.30 & 1.68 &   39 &  \textbf{800306.64} &      0.97 &            0.74 &       70 &      42 \\ 
\hline
\end{tabular}
\label{tab:lab1_3600s}
\end{table}

\end{landscape}

\subsection{Results for G3LSPD-C}
\label{sec:resultsg3lspdc}

In this section, we assess the performance of our hybrid heuristic when applied to G3LSPD-C and analyze how it behaves under different capacity configurations. We did not perform any specific parameterization for RFFO, instead, we used the same parameters which were applied for 3LSPD-C (see Section~\ref{sec:approachesandsettings}).
RFFO is compared against ES-LS, with both approaches executed with a time limit of 600 seconds for each of the instances.

The results are summarized in Tables \ref{tab:unbal_capX}-\ref{tab:bal_capX}.
The first three columns identify the capacity configurations of the instances, given by the capacity factors at the plant ($C$), at the warehouses ($C_s^{W}$), and at the retailers ($C_s^{R}$), followed by the number of instances with that specific configuration.
A value $\infty$ indicates that there is no limited capacity.
The remaining columns are similar to those in Tables~\ref{tab:lab0_600s}-\ref{tab:lab1_600s}.
The results show that RFFO remains very robust for G3LSPD-C. Besides, it clearly outperforms ES-LS when the same time limit of 600 seconds is available for both approaches.
RFFO obtained best average values for all the capacity configurations but one for the unbalanced instances ($C=\infty$,$C^W_s=\infty$,$C^R_s=1.50$) and two for the balanced ones ($C=\infty$,$C^W_s=\infty$,$C^R_s=1.50$; and $C=2.00$,$C^W_s=\infty$,$C^R_s=1.50$). Overall, a solution at least matching the best obtained by ES-LS was achieved by RFFO for around 79\% of the instances (2272 and 2270 out of 2880, for the unbalanced and balanced, respectively) and strictly better solutions were achieved for around 48\% (1374 out of 2880) of the unbalanced instances and around 56\% (1607 out of 2880) of the balanced ones.
It can also be noticed that the average optimality gaps of the solutions encountered by RFFO are reasonably low, with the largest values being 2.60\% for the unbalanced instances and 2.77\% for the balanced ones. 
It is noteworthy that RFFO achieved good average improvements over ES-LS for the great majority of the capacity configurations, reaching values as high as 10.91\% for the unbalanced instances and 12.54\% for the balanced ones.

\begin{figure}[H]
\centering
   \includegraphics[width=0.55\textwidth] {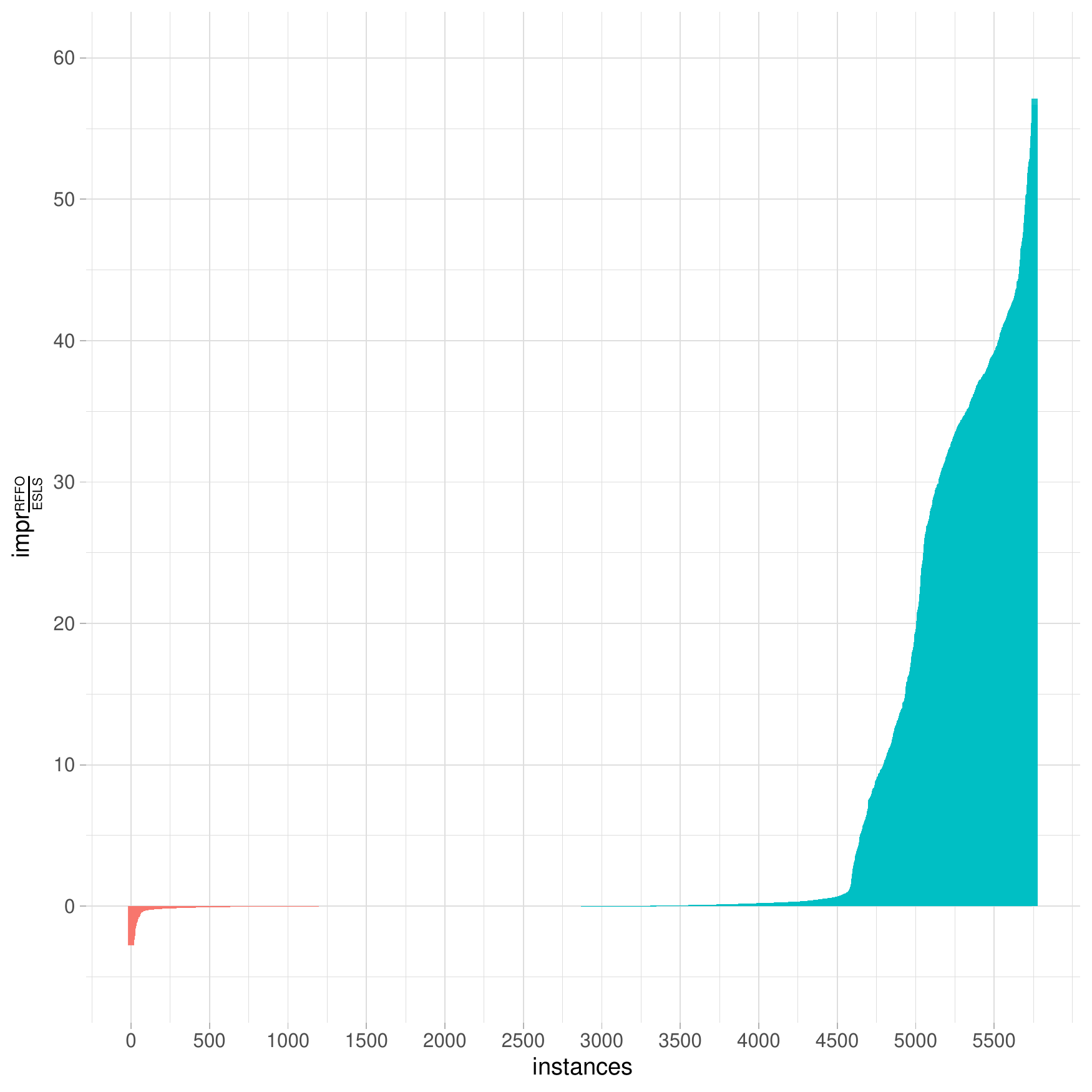}
\caption{Plot with the improvements of RFFO over ES-LS (impr$_{\frac{RFFO}{ESLS}}$) for all G3LSPD-C instances.}
\label{fig:improvementsgeneralized}
\end{figure}

The plot in Figure~\ref{fig:improvementsgeneralized} shows the improvements of RFFO over ES-LS for all the 5760 instances (both balanced and unbalanced). The plot shows that the positive improvements of RFFO over ES-LS are remarkable, as they reach values over 50\%. This contrasts with the improvements of ES-LS over RFFO, which are very modest.

\begin{landscape}

\begin{table}[H]
\centering
\footnotesize
\caption{Results for G3LSPD-C considering the unbalanced instances and ES-LS running for 600 seconds.}
\begin{tabular}{l c c c|ccc|cccccccc} \hline
\multicolumn{4}{c|}{Instances}  &  \multicolumn{3}{c|}{ES-LS} &  \multicolumn{8}{c}{RFFO} \\
 $C$ &  $C_s^{W}$ &  $C_s^{R}$ & \#inst &    best$_{ESLS}$ &  gap$_{ESLS}$ &  \#opt &   best$_{RF}$  &   best$_{RFFO}$  &  gap$_{RFFO}$ &  impr$_{\frac{FO}{RF}}$ &  impr$_{\frac{RFFO}{ESLS}}$ &  $\leq$best$_{ESLS}$ & $<$best$_{ESLS}$ &  rnd$_{FO}$ \\ \hline
  ALL &   ALL &   ALL &  2880 &  834685.45 &  6.48 &  960 &  816925.90 &  \textbf{743460.02} &      1.34 &        5.24 &            5.33 &     2272 &    1374 &    9.83 \\  
  $\infty$ &   $\infty$ &   1.50 &   120 &  \textbf{484160.51} &  0.53 &   90 &  543917.91 &  484311.72 &      0.55 &        5.01 &           -0.02 &      101 &      10 &   15.51 \\
  $\infty$ &   $\infty$ &  1.75 &   120 &  541489.87 &  4.87 &   76 &  559477.57 &  \textbf{480507.27} &      0.86 &        6.10 &            4.22 &      110 &      29 &   14.76 \\
  $\infty$ &   $\infty$ &   2.00 &   120 &  544119.62 &  5.27 &   72 &  559475.95 &  \textbf{478024.94} &      0.98 &        6.37 &            4.52 &      103 &      29 &   14.40 \\
  $\infty$ &   1.50 &   $\infty$ &   120 &  607501.86 &  8.33 &   49 &  574476.13 &  \textbf{483162.53} &      1.19 &        7.00 &            7.48 &      102 &      42 &    9.33 \\
  $\infty$ &  1.75 &   $\infty$ &   120 &  603603.63 &  8.18 &   52 &  571461.91 &  \textbf{482099.47} &      1.19 &        6.96 &            7.31 &      101 &      39 &    9.79 \\
  $\infty$ &   2.00 &   $\infty$ &   120 &  600162.89 &  7.88 &   61 &  572080.80 &  \textbf{480906.80} &      1.12 &        6.97 &            7.09 &      104 &      40 &    9.90 \\
  1.50 &   $\infty$ &   1.50 &   120 &  902756.38 &  0.67 &   34 &  943362.97 &  \textbf{901812.56} &      0.60 &        2.70 &            0.07 &       84 &      55 &   10.74 \\
  1.50 &   $\infty$ &  1.75 &   120 &  932684.83 &  2.89 &   26 &  948801.66 &  \textbf{898322.14} &      0.76 &        3.30 &            2.17 &       84 &      63 &   10.02 \\
  1.50 &   $\infty$ &   2.00 &   120 &  937679.52 &  3.49 &   23 &  951317.95 &  \textbf{896930.26} &      0.89 &        3.56 &            2.66 &       87 &      71 &    9.59 \\
  1.50 &   1.50 &   $\infty$ &   120 & 1069909.11 & 10.15 &   11 &  994423.41 &  \textbf{899323.09} &      1.10 &        6.06 &            9.24 &       84 &      73 &    6.62 \\
  1.50 &  1.75 &   $\infty$ &   120 & 1053151.71 &  9.43 &   11 &  985842.44 &  \textbf{898659.68} &      1.09 &        5.63 &            8.53 &       91 &      79 &    6.73 \\
  1.50 &   2.00 &   $\infty$ &   120 & 1042370.48 &  8.96 &   10 &  985806.15 &  \textbf{898482.92} &      1.14 &        5.54 &            8.00 &       86 &      74 &    6.90 \\
 1.75 &   $\infty$ &   1.50 &   120 &  837564.78 &  1.94 &   51 &  879762.05 &  \textbf{836332.02} &      1.84 &        3.02 &            0.10 &       93 &      49 &   11.13 \\
 1.75 &   $\infty$ &  1.75 &   120 &  854492.82 &  3.49 &   45 &  887781.96 &  \textbf{832540.36} &      1.99 &        3.81 &            1.56 &       96 &      56 &   10.76 \\
 1.75 &   $\infty$ &   2.00 &   120 &  857483.93 &  3.95 &   41 &  887351.69 &  \textbf{830461.54} &      2.06 &        3.94 &            1.97 &       97 &      64 &   10.24 \\
 1.75 &   1.50 &   $\infty$ &   120 & 1016222.15 & 12.28 &   17 &  919156.36 &  \textbf{834238.27} &      2.60 &        5.81 &           10.14 &       93 &      77 &    6.79 \\
 1.75 &  1.75 &   $\infty$ &   120 & 1006665.90 & 11.83 &   25 &  915495.80 &  \textbf{833203.91} &      2.56 &        5.68 &            9.71 &       94 &      75 &    7.13 \\
 1.75 &   2.00 &   $\infty$ &   120 &  988828.15 & 10.82 &   24 &  912322.79 &  \textbf{832210.11} &      2.51 &        5.47 &            8.69 &       97 &      78 &    7.09 \\
  2.00 &   $\infty$ &   1.50 &   120 &  764183.61 &  0.84 &   57 &  813729.07 &  \textbf{762869.21} &      0.73 &        3.68 &            0.12 &       96 &      45 &   12.02 \\
  2.00 &   $\infty$ &  1.75 &   120 &  777525.17 &  2.30 &   49 &  822126.56 &  \textbf{759860.60} &      0.96 &        4.46 &            1.38 &       97 &      56 &   11.71 \\
  2.00 &   $\infty$ &   2.00 &   120 &  778023.04 &  2.59 &   46 &  823774.00 &  \textbf{757424.64} &      1.00 &        4.77 &            1.64 &       95 &      58 &   11.09 \\
  2.00 &   1.50 &   $\infty$ &   120 &  952165.01 & 12.05 &   26 &  854816.90 &  \textbf{761682.45} &      1.56 &        6.74 &           10.86 &       95 &      74 &    7.65 \\
  2.00 &  1.75 &   $\infty$ &   120 &  952370.36 & 12.02 &   32 &  850168.55 &  \textbf{760366.15} &      1.45 &        6.58 &           10.91 &       88 &      68 &    7.82 \\
  2.00 &   2.00 &   $\infty$ &   120 &  927335.35 & 10.78 &   32 &  849291.09 &  \textbf{759307.93} &      1.41 &        6.54 &            9.66 &       94 &      70 &    8.08 \\ \hline
\end{tabular}
\label{tab:unbal_capX}
\end{table}

\begin{table}[H]
\centering
\footnotesize
\caption{Results for G3LSPD-C considering the balanced instances and ES-LS running for 600 seconds.}
\begin{tabular}{l c c c|ccc|cccccccc} \hline
\multicolumn{4}{c|}{Instances}  &  \multicolumn{3}{c|}{ES-LS} &  \multicolumn{8}{c}{RFFO} \\
 $C$ &  $C_s^{w}$ &  $C_s^{r}$ & \#inst &    best$_{ESLS}$ &  gap$_{ESLS}$ &  \#opt &   best$_{RF}$  &   best$_{RFFO}$  &  gap$_{RFFO}$ &  impr$_{\frac{FO}{RF}}$ &  impr$_{\frac{RFFO}{ESLS}}$ &  $\leq$best$_{ESLS}$ & $<$best$_{ESLS}$ &  rnd$_{FO}$ \\ \hline
ALL &   ALL &   ALL &  2880 &  890596.63 &  6.71 &  600 &  862395.36 &  \textbf{785533.74} &      1.30 &        5.25 &            5.56 &     2270 &    1607 &    9.44 \\
  $\infty$ &   $\infty$ &   1.50 &   120 &  \textbf{508796.98} &  0.09 &   96 &  568798.50 &  508932.72 &      0.10 &        4.57 &           -0.01 &      108 &       3 &   16.32 \\
  $\infty$ &   $\infty$ &  1.75 &   120 &  508272.09 &  0.58 &   90 &  578266.18 &  \textbf{504945.44} &      0.22 &        5.51 &            0.38 &      108 &      11 &   16.14 \\
  $\infty$ &   $\infty$ &   2.00 &   120 &  529881.64 &  2.49 &   84 &  580879.10 &  \textbf{502518.49} &      0.31 &        5.83 &            2.22 &      112 &      20 &   15.84 \\
  $\infty$ &   1.50 &   $\infty$ &   120 &  653862.10 &  8.30 &   23 &  608459.06 &  \textbf{507482.72} &      0.71 &        7.24 &            7.77 &       97 &      64 &    7.52 \\
  $\infty$ &  1.75 &   $\infty$ &   120 &  648993.73 &  8.28 &   26 &  600896.58 &  \textbf{506415.12} &      0.65 &        7.04 &            7.87 &       98 &      63 &    7.80 \\
  $\infty$ &   2.00 &   $\infty$ &   120 &  632063.21 &  7.75 &   31 &  597312.06 &  \textbf{503821.95} &      0.57 &        6.88 &            7.33 &      102 &      59 &    8.33 \\
  1.50 &   $\infty$ &   1.50 &   120 &  954138.85 &  0.75 &   14 &  994448.27 &  \textbf{953085.25} &      0.67 &        2.55 &            0.08 &       94 &      75 &    9.56 \\
  1.50 &   $\infty$ &  1.75 &   120 &  970392.93 &  2.12 &    8 & 1009293.15 &  \textbf{949128.89} &      0.77 &        3.63 &            1.38 &       85 &      71 &    9.03 \\
  1.50 &   $\infty$ &   2.00 &   120 &  974329.76 &  2.56 &    4 & 1008569.90 &  \textbf{946729.74} &      0.83 &        3.73 &            1.77 &       86 &      79 &    8.70 \\
  1.50 &   1.50 &   $\infty$ &   120 & 1159592.63 & 11.40 &    2 & 1053885.00 &  \textbf{944823.13} &      1.08 &        6.69 &           10.51 &       98 &      94 &    6.21 \\
  1.50 &  1.75 &   $\infty$ &   120 & 1153199.16 & 11.22 &    3 & 1044668.72 &  \textbf{943657.85} &      1.02 &        6.25 &           10.37 &       92 &      88 &    6.58 \\
  1.50 &   2.00 &   $\infty$ &   120 & 1134666.22 & 10.43 &    3 & 1042764.87 &  \textbf{943061.71} &      1.00 &        6.22 &            9.59 &       95 &      91 &    6.82 \\
 1.75 &   $\infty$ &   1.50 &   120 &  888180.59 &  2.07 &   30 &  927503.39 &  \textbf{887447.47} &      2.02 &        2.66 &            0.06 &       84 &      56 &   10.72 \\
 1.75 &   $\infty$ &  1.75 &   120 &  900740.97 &  3.34 &   24 &  940859.12 &  \textbf{883700.92} &      2.18 &        3.77 &            1.21 &       89 &      65 &   10.11 \\
 1.75 &   $\infty$ &   2.00 &   120 &  903605.38 &  3.78 &   19 &  943598.46 &  \textbf{881140.81} &      2.28 &        4.01 &            1.56 &       86 &      71 &   10.06 \\
 1.75 &   1.50 &   $\infty$ &   120 & 1108148.59 & 13.85 &    5 &  969422.16 &  \textbf{878998.97} &      2.77 &        6.00 &           11.55 &       98 &      89 &    6.41 \\
 1.75 &  1.75 &   $\infty$ &   120 & 1101883.31 & 13.69 &    5 &  960839.06 &  \textbf{878069.47} &      2.77 &        5.50 &           11.37 &       95 &      86 &    6.71 \\
 1.75 &   2.00 &   $\infty$ &   120 & 1079585.32 & 12.69 &    6 &  960293.46 &  \textbf{877398.22} &      2.67 &        5.57 &           10.44 &       92 &      85 &    6.98 \\
  2.00 &   $\infty$ &   1.50 &   120 &  \textbf{813229.40} &  0.93 &   38 &  861904.13 &  813603.17 &      0.95 &        3.34 &           -0.02 &       78 &      42 &   11.82 \\
  2.00 &   $\infty$ &  1.75 &   120 &  826247.29 &  2.28 &   31 &  872645.82 &  \textbf{809700.31} &      1.06 &        4.29 &            1.26 &       98 &      68 &   11.40 \\
  2.00 &   $\infty$ &   2.00 &   120 &  830062.05 &  2.85 &   26 &  874807.48 &  \textbf{807526.02} &      1.23 &        4.61 &            1.67 &       94 &      73 &   11.09 \\
  2.00 &   1.50 &   $\infty$ &   120 & 1053283.81 & 14.12 &    6 &  903104.61 &  \textbf{806453.76} &      1.92 &        6.90 &           12.54 &       95 &      90 &    7.08 \\
  2.00 &  1.75 &   $\infty$ &   120 & 1034947.27 & 13.38 &   11 &  900283.39 &  \textbf{805274.81} &      1.82 &        6.72 &           11.88 &       96 &      85 &    7.47 \\
  2.00 &   2.00 &   $\infty$ &   120 & 1006215.81 & 12.05 &   15 &  893986.26 &  \textbf{806743.57} &      1.72 &        6.47 &           10.69 &       90 &      79 &    7.90 \\
  \hline

\end{tabular}

\label{tab:bal_capX}
\end{table}

\end{landscape}

\subsection{Economic impacts of the capacity configurations}

In what follows, consider the deviation of a solution for a given capacity configuration from a baseline problem calculated as $100 \frac{\textrm{best}_{RFFO} - \textrm{best}_{\textrm{base}}}{\textrm{best}_{\textrm{base}}}$ for a specific instance. In this formula, $\textrm{best}_{RFFO}$ defines the best solution achieved by RFFO for the given capacity configuration while $\textrm{best}_{\textrm{base}}$ denotes the optimal solution for the baseline problem 3LSPD-U (determined with the approaches described in \citeA{CunMel21}) or the best solution achieved by RFFO for the baseline problem 3LSPD-C. Note that the deviation indicates the economical impact (in terms of increase in the objective value) of the given capacity configuration on the obtained solution. 

\begin{figure}[H]
\centering
   \includegraphics[width=0.5\textwidth] {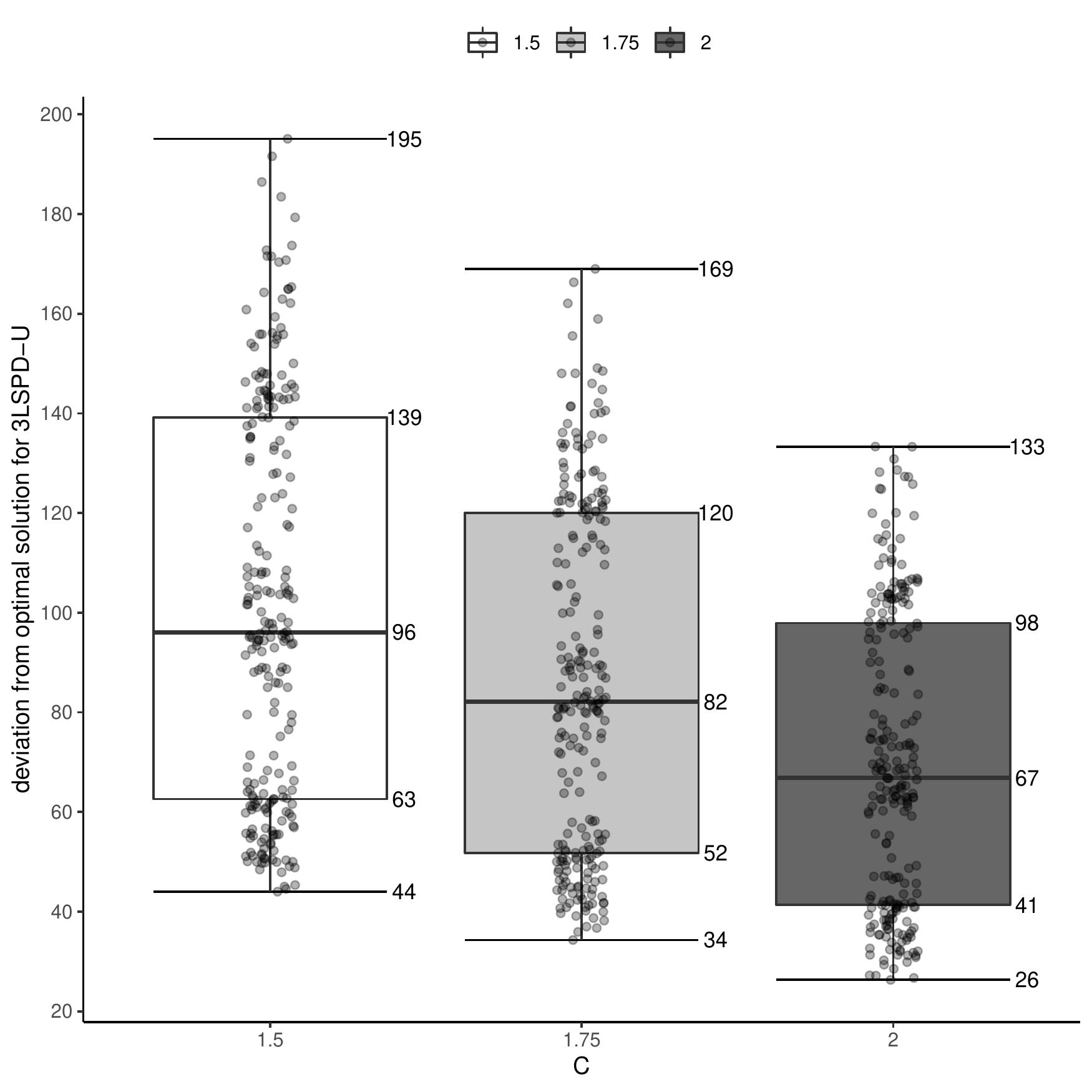}
\caption{Boxplot summarizing the deviations from the optimal solutions for the 3LSPD-U when capacities occur at the plant.}
\label{fig:boxplotuncapplantcapplant}
\end{figure}

The boxplots in Figure~\ref{fig:boxplotuncapplantcapplant} depict the deviations from the baseline problem 3LSPD-U when there are production capacities at the plant.
It can be seen that median deviations of 96\%, 82\%, and 67\% were incurred when considering production plant capacity factors of 1.50, 1.75, and 2.00, respectively. It is noteworthy that the economical impacts can reach values as high as 195\% for the tighter capacity factor.
The large deviations can be explained by the need for additional setups, which have high costs in the considered instances, when the amount that can be produced in a given period is limited. 

\begin{figure}[H]
\subfigure[Baseline problem 3LSPD-U.]{
   \includegraphics[scale =0.47] {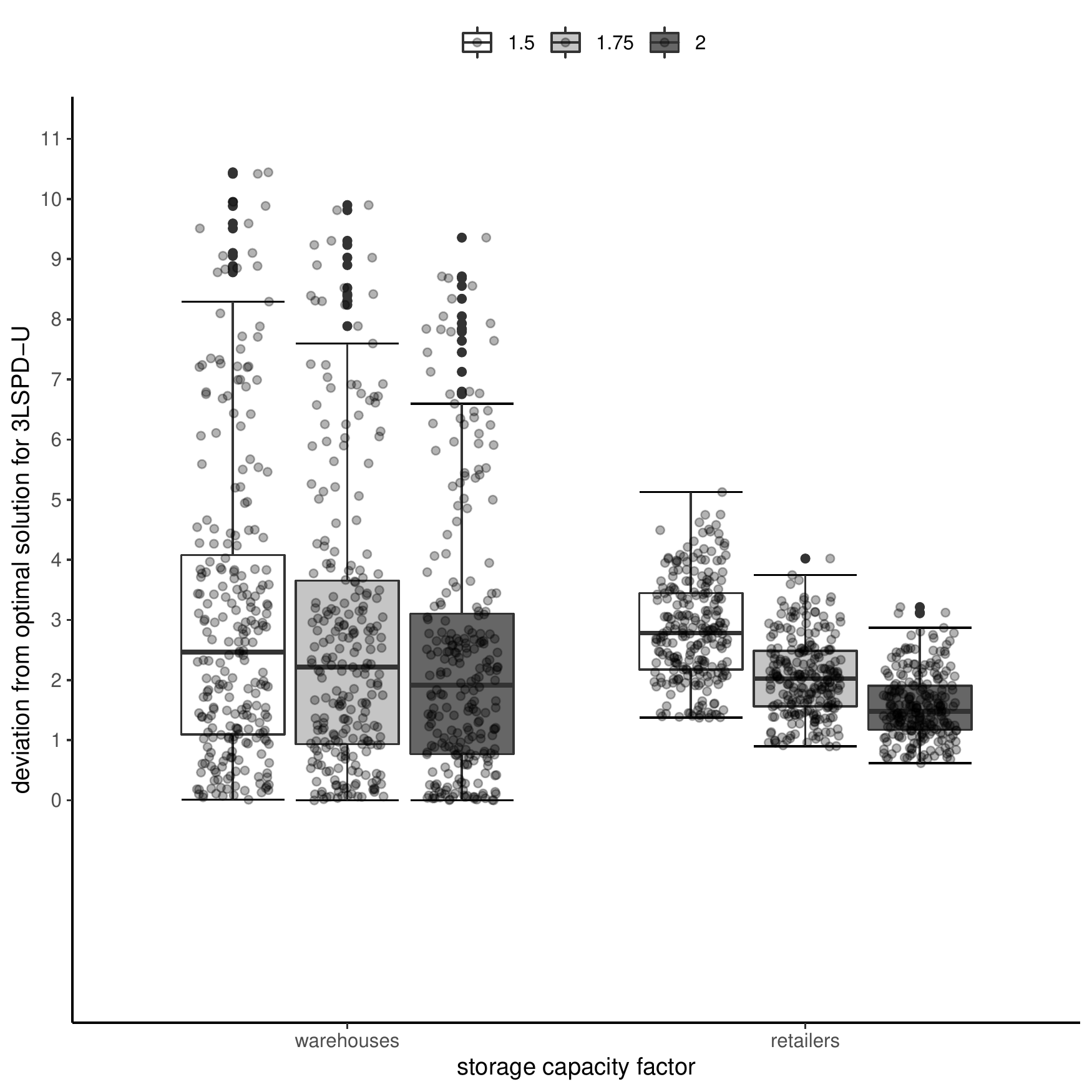}
   \label{fig:boxplotuncapplant}
 }
  \subfigure[Baseline problem 3LSPD-C with $C = 1.5$.]{
   \includegraphics[scale =0.47] {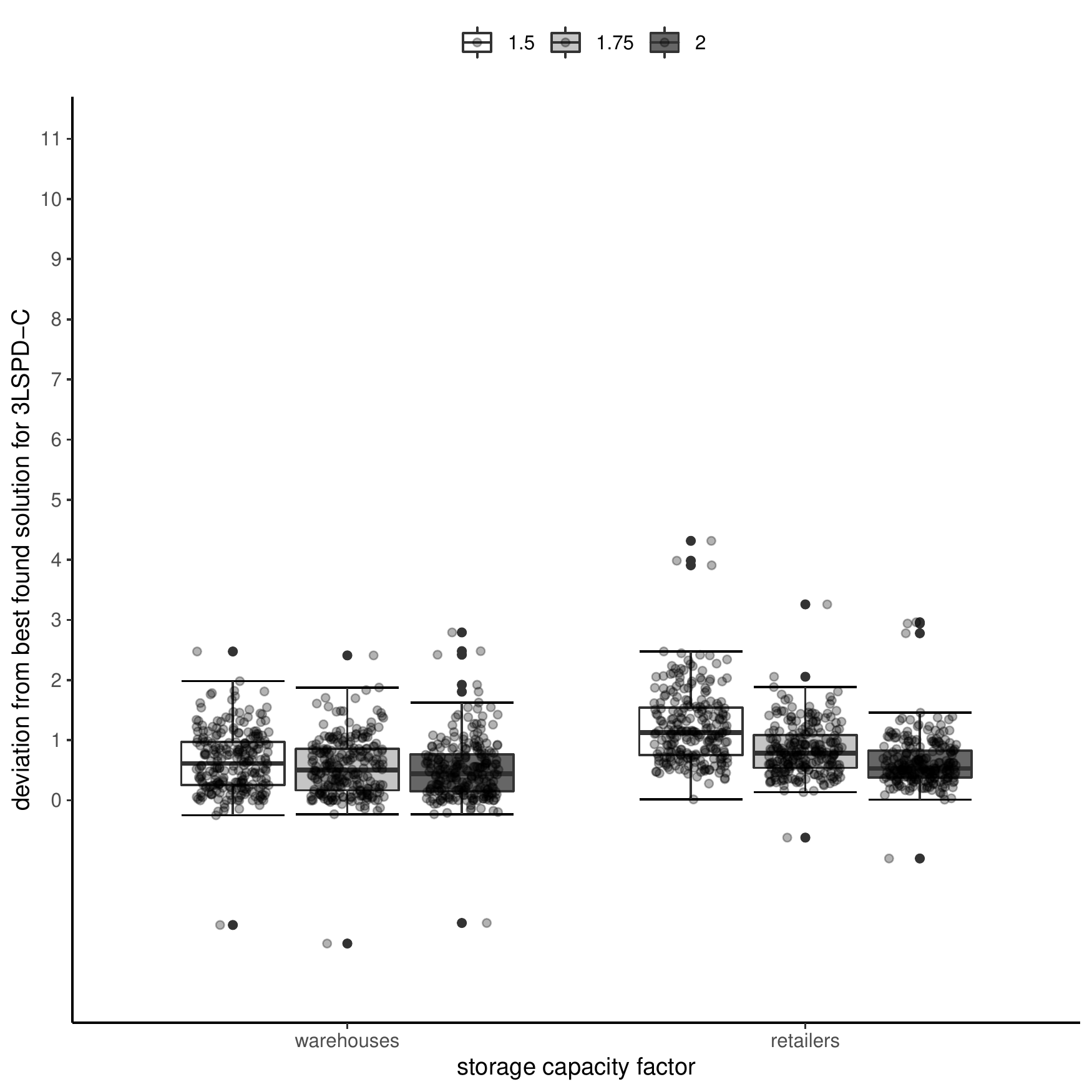}
   \label{fig:boxplotplant1_5}
 }
  \subfigure[Baseline problem 3LSPD-C with $C = 1.75$.]{
   \includegraphics[scale =0.47] {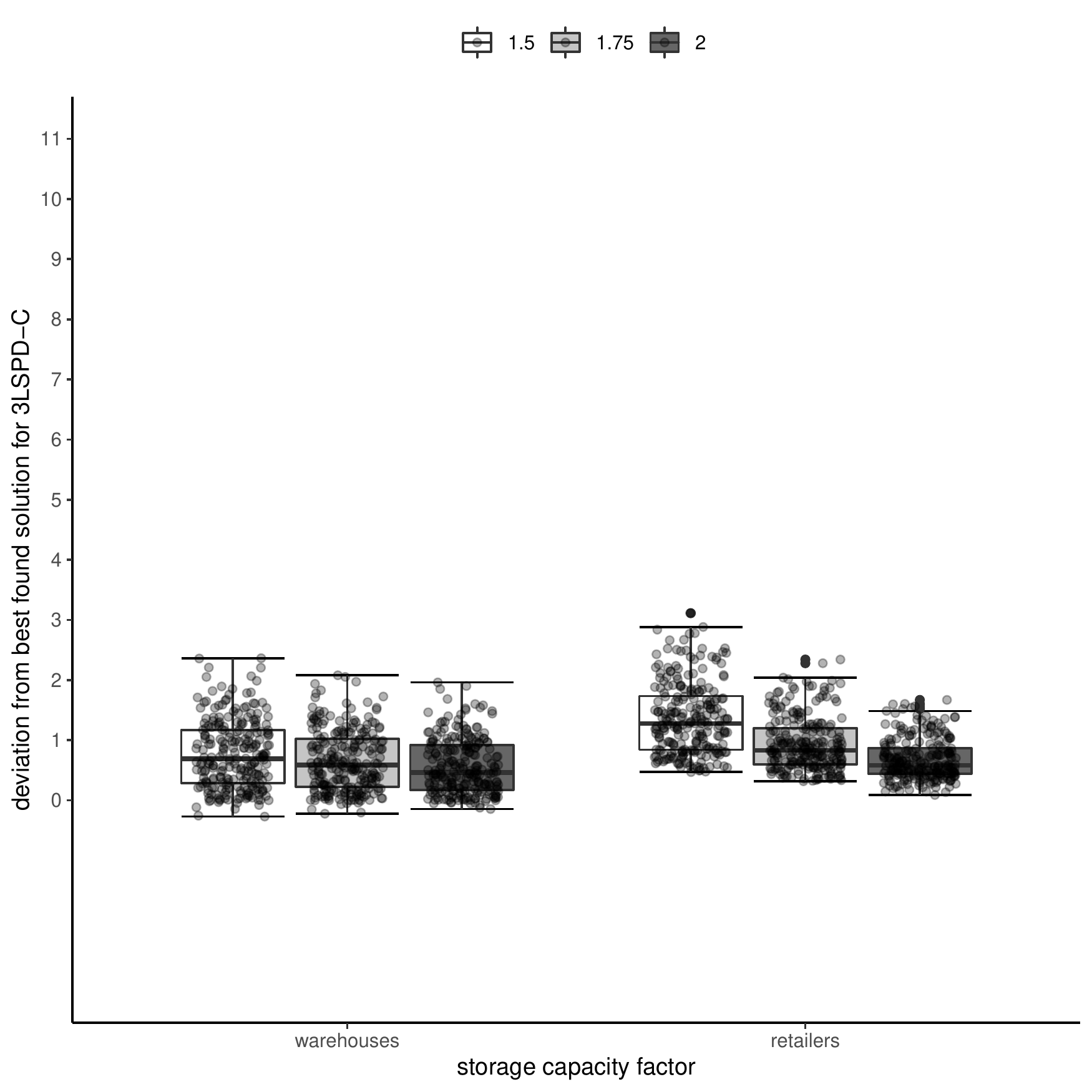}
   \label{fig:boxplotplant1_75}
 }
  \subfigure[Baseline problem 3LSPD-C with $C = 2.00$.]{
   \includegraphics[scale =0.47] {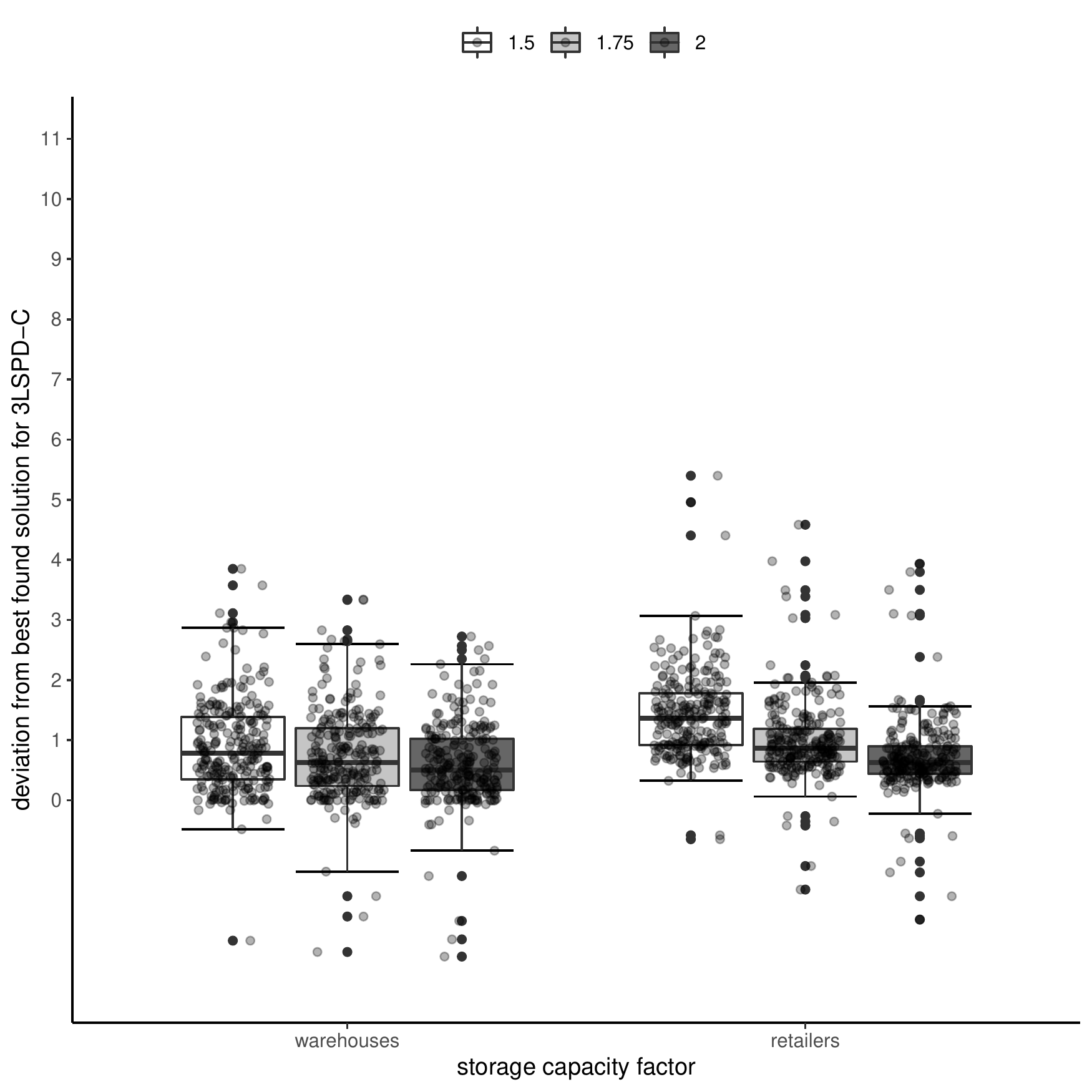}
   \label{fig:boxplotplant2}
 }
\caption{Boxplot summarizing the deviations from the best known solutions for the baseline problems when capacities occur at the warehouses or at the retailers.}
\label{fig:boxplotuncapplantcapwr}
\end{figure}

The boxplots in Figure~\ref{fig:boxplotuncapplantcapwr} show the deviations from the baseline problems 3LSPD-U (Subfigure~\ref{fig:boxplotuncapplant}) and 3LSPD-C (Subfigures~\ref{fig:boxplotplant1_5}-\ref{fig:boxplotplant2}) when there are limited storage capacities at the warehouses or at the retailers.
Notice that all subfigures use the same scale.
Overall, the boxplot in Subfigure~\ref{fig:boxplotuncapplant} shows that when there is no production capacity at the plant, there is a larger variance on the deviations for the capacity factors at the warehouses, implying reasonably high deviations for several instances. On the other hand, when the baseline problem is 3LSPD-C, the deviations are more modest and stable, with the capacity factors at the retailers implying solutions with higher costs.
Notice that some negative values occur in Subfigures~\ref{fig:boxplotplant1_5}-\ref{fig:boxplotplant2}, which can be explained by the fact that the hybrid heuristic RFFO can in certain situations tackle some more restricted subproblems of G3LSPD-C better than those of 3LSPD-C, achieving solutions that are closer to optimal.

\section{Concluding remarks}
\label{sec:concludingremarks}

In this paper, we considered the capacitated three-level lot-sizing and replenishment problem with a distribution structure (3LSPD-C) and an extension denoted the generalized capacitated three-level lot-sizing and replenishment problem with a distribution structure (G3LSPD-C) in which storage capacities restrict the operations on the warehouses and/or retailers. We proposed a very effective hybrid MIP heuristic combining relax-and-fix and fix-and-optimize based on varying-size neighborhoods (RFFO) that can be applied for both problems. 

Extensive computational experiments showed that RFFO clearly outperforms a state-of-the-art formulation (ES-LS) in terms of solutions quality when both approaches are offered the same time limit.
Depending on the instances' characteristics, average improvements can reach approximately 8\% for 3LSPD-C.
For G3LSPD-C, depending on the capacity configurations, average improvements can come to more than 12\%.
Such observed improvements can allow considerable economical impacts in a practical setting.
The results also show that the solutions encountered by RFFO can be very competitive even when compared to ES-LS running for a much larger time.
Overall, RFFO is very robust in that it obtains, for both 3LSPD-C and G3LSPD-C, solutions whose average optimality gaps are below 3\% for all the reported instance configurations.

We remark that our approach can be easily extended to other variants and extensions of the problem, such as those involving multiple plants, backlogging, retailers that can receive items from multiple warehouses, among others. Furthermore, we believe that the concepts used in our fix-and-optimize approach may bring enhancements when applied to other problems involving supply chains and/or multi-period logistical optimization problems.



\section*{Acknowledgments} 
Work of Rafael A. Melo was supported by the State of Bahia Research Foundation (FAPESB) and the Brazilian National Council for Scientific and Technological Development (CNPq).
Work of Geraldo R. Mateus was supported by Fundação de Amparo à Pesquisa do Estado de Minas Gerais (FAPEMIG) and CNPq.

\bibliographystyle{apacite}

\bibliography{00_main}

\newpage

\begin{appendices}

\section{Echelon stock formulation}
\label{sec:echelonstockformulations}

In this section, we present the echelon stock formulation that is employed in our hybrid heuristic. This formulation was introduced in \citeA{GruBazCorJan19}, where more details can be obtained.
Define the echelon stock $I^i_t$ for facility $i\in F$ in period $t\in T$ as
\begin{equation*}
    I^i_t = 
    \begin{cases}
     s^i_t + \sum_{w\in W} s^w_t + \sum_{r\in R} s^r_t,  & \textrm{if } i=p;\\
     s^i_t + \sum_{r \in \delta(i)}, & \textrm{if } i\in W;\\
     s^i_t,  & \textrm{if } i\in R.
    \end{cases}
\end{equation*}
The formulation ES-LS can be cast as
\begin{align}
z_{ES\textrm{-}LS} = & \  \min \ \ \ \sum_{t \in T} \left(  \sum_{i\in F} sc^i_t y^i_t + hc^p_t I^p_t +  \sum_{w\in W} (hc^w_t - hc^p_t ) I^w_t + \sum_{ r \in R} (hc^r_t - hc^{\delta_w(r)}_t ) I^r_t \right)  \label{es-obj} & \\
(ES\textrm{-}LS) \qquad & \eqref{std-3}\textrm{-}\eqref{std-3b}, \eqref{std-5}, \nonumber \\
& I^i_{t-1} + x^i_t = d^i_t + I^i_t, \qquad  \textrm{for} \ i \in F, \ t\in T, \label{es-1} \\
&   I^i_t \geq \sum_{j \in \delta(i)} I^j_t, \qquad  \textrm{for} \ i \in P \cup W, \ t\in T, \label{es-ww}\\
&   I^i_{t-1} \geq \sum_{j = t}^l d^i_j \left( 1 - \sum_{u=t}^j y^i_u \right), \qquad  \textrm{for} \ i \in F, \ t,l\in T, \textrm{ with }t \leq l , \label{es-2}\\
&  x^i_{t}, \ I^i_t \geq 0, \qquad  \textrm{for} \ i \in F, \ t\in T. \label{es-3}
\end{align}

\end{appendices}
 
\end{document}